\documentclass[11pt,leqno,oneside,letterpaper]{amsart}

\usepackage{amssymb,amsmath,amsfonts,latexsym, amscd, amsthm, url} 
\usepackage[letterpaper,left=1.75truein, top=1truein]{geometry}
\usepackage{color}
\usepackage[colorlinks,linkcolor=blue,citecolor=blue]{hyperref}

\usepackage{caption,subcaption,graphicx,epsfig}
\usepackage{tikz-cd,calc}
\usepackage{epsfig}
\usepackage{graphicx}
\usepackage{tikz-cd}
\usepackage{array}
\usepackage{float}
\usetikzlibrary{knots}
\usetikzlibrary{decorations.markings}
\definecolor{light-gray}{gray}{0.92}
\definecolor{ultra-light-gray}{gray}{0.97}

\graphicspath{{images/}}

\newtheorem{theorem}{Theorem}[section]
\newtheorem{lemma}[theorem]{Lemma}

\newtheorem{definition}[theorem]{Definition}

\newtheorem{corollary}[theorem]{Corollary} 
\newtheorem{claim}[theorem]{Claim}
\theoremstyle{definition}
\newtheorem{remark}[theorem]{Remark}

\newtheorem{conjecture}[theorem]{Conjecture}

\newtheoremstyle{cases}
  {12pt plus 6 pt}
  {2pt}
  {\bfseries}   
  {}
  {\bfseries}
  {.}
  {.5em}
  {}
\theoremstyle{cases}

\numberwithin{subcase}{case} \numberwithin{subsubcase}{subcase}
\numberwithin{equation}{subsection}

\address{D\'epartement de Math\'ematiques, Universit\'e du Qu\'ebec \`a Montr\'eal, 201 Avenue du Pr\'esident-Kennedy, 
Montr\'eal,
QC H2X 3Y7.} 
\email{ba.idrissa@courrier.uqam.ca}
 
\title{L-spaces, left-orderability and two-bridge knots}
\author{Idrissa Ba}
\begin{document}

\maketitle

\begin{center}
\today
\end{center}

\begin{abstract}We show that the 3-fold cyclic branched cover of any genus 2 two-bridge knot
$K_{[-2q,2s,-2t,2l]}$ is an L-space and its fundamental group is not left-orderable. Therefore 
the family of 3-fold cyclic branched cover of any genus 2 two-bridge
knot $K_{[-2q,2s,-2t,2l]}$ verifies the $L$-space conjecture.
We also show that if  $K_{[2k,-2l]}$ is a 2-bridge knot with $k\geq 2$, $l>0$,
then the fundamental group of the 5-fold cyclic branched cover of  $K_{[2k,-2l]}$ is not left-orderable,
which will complete the proof that the fundamental group of the 5-fold cyclic branched cover of any genus one two-bridge knot
is not left-orderable.
\end{abstract}

\section{Introduction}
In this paper we study the L-space conjecture for the cyclic branched covers of low genus two-bridge knots.

A closed, connected 3-manifold $M$ is an {\it L-space} if it is a rational homolgy sphere 
with the property that $rk\widehat{HF}(M)=ord(H_1(M,\mathbb{Z}))$
 ([OSz], [OSz06]).
 
A group $G$ is called {\it left-orderable} if there exists a strict total ordering $<$, of $G$ such that $g<h$ implies $fg<fh$ for 
all $f$, $g$, $h$ $\in$ $G$. By convention the trivial group is not left-orderable.

A closed, connected, orientable 3-manifold $M$ is called a {\it total $L$-space} if it is an $L$-space whose fundamental group is 
not left-orderable.
\begin{conjecture}{\rm (Conjecture 1 in [BGW])}
An irreducible rational homolgy 3-sphere is an L-space if and only if
 its fundamental
group is not left-orderable.
\end{conjecture}
Let $K_{[a_1, a_2, \cdots, a_m]}$ denote the two-bridge knot of type $\frac{p}{q}$, where $[a_1, a_2, \cdots, a_m]$ is a 
continued fraction expansion for $\frac{p}{q}$. We follow the convention that 
$\frac{p}{q}=a_1+\frac{1}{a_2+\cdots +\frac{1}{a_m}}$. Every two-bridge 
knot admits a continued fraction expansion with an even number of even parameters $[2a_1, 2b_1, 2a_2,2 b_2,  \cdots, 2a_m, 2b_m]$.
The two-bridge knots $K_{[2a_1, 2b_1, 2a_2,2 b_2,  \cdots, 2a_m, 2b_m]}$ are of genus $m$, when the 
$a_i$ and $b_i$ are in $\mathbb{Z}\setminus\{0\}$. 
Every genus 2 two-bridge knot can be written as $K_{[2a_1, 2b_1, 2a_2,2 b_2]}$, where
$a_i$ and $b_i$ are in $\mathbb{Z}\setminus\{0\}$, $i=1,2$.

If $K$ is a genus one, alternating knot, then $K$ is either a genus one two bridge knot or, up to mirroring, a pretzel knot
$P(2n+1, 2m+1, 2p+1)$ with $m$, $n$, $p$ positive integers ([BZ], Lemma 3.1). In the case where $K$ is a 
genus one two bridge knot, much work has been done to study the left-orderability of the fundamental groups of
its cyclic branched covers $\pi_1(\Sigma_n(K))$. In this direction, the fundamental group of the 2-fold branched cover of
$K_{[2k,-2l]}$ is not left-orderable for $k>0$
and $l>0$, because $\Sigma_2(K_{[2k,-2l]})$ is a lens space, and in [DPT] it is shown that the fundamental group of
the 3-fold cyclic branched cover of  $K_{[2k,-2l]}$ 
is not left-orderable for $k>0$
and $l>0$. Gordon and Lidman [GL] showed that the fundamental group of the 4-fold cyclic branched cover 
of  $K_{[2k,-2l]}$ is not left-orderable for $k>0$
and $l>0$. But this is false for $n$ sufficiently large by [Hu] and [Tra].

The knot $5_1$ corresponds to the two-bridge knot $K_{[-2,2,-2,2]}$. 
Since the 3-fold cyclic branched cover of $5_1$ is the Poincar\'e homology 
sphere, it is a total $L$-space. Therefore
we can ask the 
following question: Is the 3-fold cyclic branched cover of the knots $K_{[-2q,2s,-2t,2l]}$ a total L-space?
In this paper we answer this question positively.

\begin{theorem}\label{thm: main result4}
  The 3-fold cyclic branched cover of $K_{[-2q,2s,-2t,2l]}$ is a total L-space, where $q$, $s$, $t$ and $l$ 
  $\in\mathbb{Z}\setminus\{0\}$.
\end{theorem}
\begin{theorem} \label{thm: main result3}
For $k\geq 2$, $l>0$, the fundamental group of
 the 5-fold cyclic branched cover of  $K_{[2k,-2l]}$ is not left-orderable.
\end{theorem}

Theorem \ref{thm: main result3} combines with Theorem 2 in [DPT] and the result of Mitsunori Hori see [Te],
to imply the following corollary.
\begin{corollary}
The 5-fold cyclic branched cover of any genus one two-bridge knot is a total L-space.
\end{corollary}

The paper is organized as follows. In the second section we introduce some background material and notations. In section 3 we 
proof Theorem \ref{thm: main result3}. In section 4, we prove that
the fundamental group of the 3-fold cyclic branched cover of any genus 2 two-bridge knot is not left-orderable.
Finaly, in section 5 we prove that the 3-fold cyclic branched cover of any genus 2 two-bridge knot is an L-space, thus
completing the proof of Theorem \ref{thm: main result4}.

{\bf{Acknowledgment.}} I would like to thank my supervisor Professor Steven Boyer for drawing my attention
to the topic of the current paper and his consistent encouragement and support.
\section{Background notions, terminology, and notation}

In this section we define some basic notions which will be useful in this paper.

Let $K$ be an oriented knot in $\mathbb{S}^3$. Let $M_K$ be the exterior of $K$ and $S$ be a Seifert surface for $K$.
Isotope $S$ so that $S\cap\partial M_K$ is a longitude of $K$ and
let $F=S\cap M_K$. Let $C$ be a tubular neighborhood of $F$ in $M_K$.
Then $C$ is homeomorphic to $F\times [-1,1]$. Let $Y:=M_K-int(C)$. The boundary of $Y$ has two copies
${F}^{-}\cong F\times \{-1\}$ and ${F}^+\cong F\times \{1\}$. We have a triple $(Y,F^+,F^-)$. Consider $n$-copies of 
this, denoted by $(Y_i,F_i^+,F_i^-)$, $i=0,\cdots, n-1$, and glue them together by identifying $F_0^+\subset Y_0$ with
 $F_1^-\subset Y_1$, $F_1^+\subset Y_1$ with $F_2^-\subset Y_2$, $\cdots$,  $F_{n-2}^+\subset Y_{n-2}$ with  
 $F_{n-1}^-\subset Y_{n-1}$ and $F_{n-1}^+\subset Y_{n-1}$ with $F_0^-\subset Y_0$. Call the resulting space $Y_n$. 
There is a regular covering map $g:Y_n\longrightarrow M_K$ and its 
 group of deck transformations is isomorphic to $\mathbb{Z}_n$. 
The manifold $Y_n$ is called the $n$-fold cyclic cover of $M_K$ and its 
 fundamental group isomorphic to $Ker(\pi_1(M_K)\longrightarrow \mathbb{Z}_n)$. 
To construct the $n$-fold cyclic branched
 cover $\Sigma_n(K)$, we have to glue a solid torus $V\cong D^2\times \mathbb{S}^1$ to $Y_n$ by identifying the meridian
 $\partial D^2\times \{1\}$ of $V$ with the preimage of the meridian $\mu$ of $\partial M_K$  under $g:Y_n\longrightarrow M_K$.
The manifold $\Sigma_n(K)$ is a closed oriented 3-manifold.

For the construction of the $n$-fold cyclic branched cover of an oriented link $L$ see [BBG].
\begin{definition}{\rm
Let L be a link and $D$ a link diagram of $L$. Checkerboard color the regions of the 
complement of the diagram in $\mathbb{R}^2$. Assume that the unbounded region 
$X_0$ is
colored white. The other white regions will be called by $X_1$, $X_2$, $\cdots$, $X_n$. To any crossing $p$ of L we 
associate the number $\chi(p)$ which is +1 or -1 according to the convention in the Figure 1.
\begin{figure}[H]
\centering
\def\svgwidth{0.30\columnwidth}
 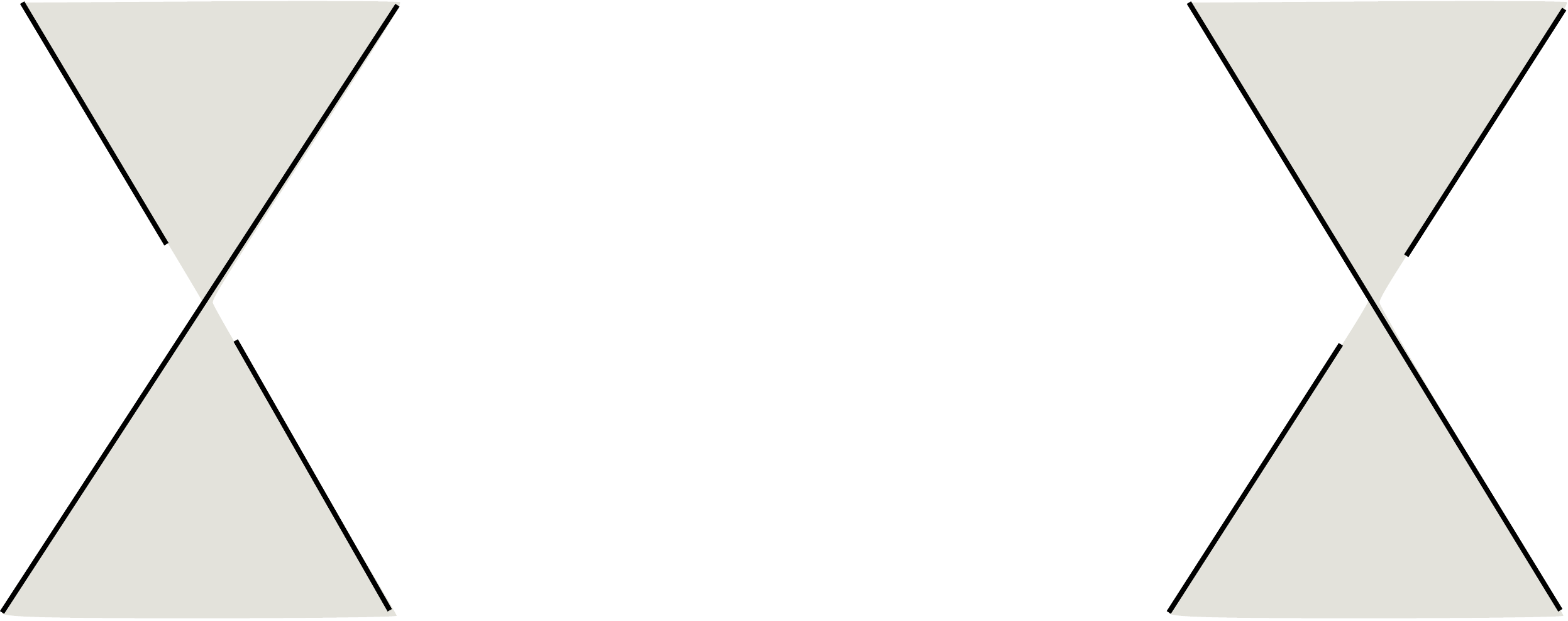
 \caption{The signs convention +1 and -1 respectively.}
\label{fig:The link A.}
\end{figure}

Let $H=(h_{ij})_{i,j=0, 1,\cdots, n}$, where

$h_{ij}$ $ = \begin{cases} -\sum_p\chi(p)$, if $i\neq j$ and the summation extends over all crossings which connect $X_i$ and
$X_j$ 
$\\-\sum_{k=0; k\neq i}^n  h_{ik}$        if $i=j$  $\end{cases}$ 

The matrix $H$ is called the {\it unreduced Goeritz matrix} of $D$. 
The {\it Goeritz matrix} $G$ of $D$ is obtained from $H$ by removing the 
first row and the first column of $H$.}
\end{definition}

Recall that the {\it determinant} of a link $L$ is the order of the first homology of its 2-fold branched cover.

\begin{theorem} Let $L$ be a non-split link. The determinant of $L$ is given by 
 the absolute value of the determinant of $G$ {\rm($\mid$det $G|$)}.
\end{theorem}

\begin{figure}[H]
\centering
\def\svgwidth{0.30\columnwidth}
 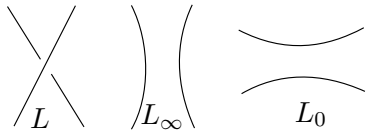
 \caption{The resolutions $L_{\infty}$ and $L_0$ respectively.}
\label{fig:The link A.}
\end{figure}

\begin{definition}{\rm
 The set $\mathcal{Q}$ of {\it quasi-alternating} links is the smallest set of links which satisfies the following properties:
 \begin{enumerate}
 \item the unknot is in $\mathcal{Q}$, 
 \item the set $\mathcal{Q}$ is closed under the following operation. Suppose $L$ is any link which admits a projection with 
 a crossing with the following properties:
 \begin{itemize}
 \item both resolution $L_0$ and $L_\infty$ are in $\mathcal{Q}$, 
 \item det$(L)$ = det$(L_0)$ + det$(L_\infty)$
\end{itemize}

\end{enumerate}
then $L$ is in $\mathcal{Q}$.}
\end{definition}

\begin{theorem} {\rm (Proposition 3.3 in [OSz])}\label{thm: 2.4}
 The 2-fold branched cover of a quasi-alternating link $L$ is an L-space.
\end{theorem}

Let $L$ be a link, $L_0$ and $L_\infty$ be the links obtained by applying the resolutions as in Figure 2 such that 
det$(L_0)\neq 0$, det$(L_\infty)\neq 0$, and det$(L)$ = det$(L_0)$ + det$(L_\infty)$. Then if 
$\Sigma_2(L_0)$ and $\Sigma_2(L_\infty)$ are $L$-spaces, then $\Sigma_2(L)$ is also an $L$-space ([OSz]).

\section{Left-orderability and genus one two-bridge knots}
In this section we will prove Theorem \ref{thm: main result3}. 
The following Lemma will be important in the proof.
\begin{lemma} \label{lemma: stabilizer}
 Let $G$ be a left-orderable group which acts by order-preserving automorphisms on a totally ordered set $(X, <_X)$.
 Then for each $a\in X$, there exists a left-order $<$ on $G$ such that the stabilizer of $a$,
 $Stab_G(a)$ is $<$-convex. The order $<$ is defined by:\\
 $g\in P(<)$ $\Leftrightarrow \begin{cases} g(a)>a$, $or \\ g(a)=a$  $ {and} $   $1<_G g  \end{cases}$

\end{lemma}
\begin{remark}
{\rm A consequence of this definition is that if $g\geq h$, then $g(a)\geq h(a)$.}
\end{remark}

\begin{proof}[Proof of Lemma \ref{lemma: stabilizer}] Routine.

\end{proof}

 Let $G=\pi_1(\Sigma_5(K_{[2k,-2l]}))$ for $k\geq 2$, $l>0$.
 We first note that $[2k,-2l]=2k+ \frac{1}{-2l}=\frac{4kl-1}{2l}=[2k-1,1,2l-1]$. By ([DPT], Proposition 2) we have that 
 $$\pi_1(\Sigma_n(K_{[2k+1,1,2l+1]}))=\{x_1,\cdots,x_n\;|\; r_1=1,\cdots, r_n=1, x_1x_2\cdots x_n=1\},$$ where
 $k>0, l>0, r_i=x_i^{-1}(x_i^{-k}x_{i+1}^{k+1}x_i^{-1})^lx_i^{-k}x_{i+1}^{k+1}((x_{i+1}^{-k}x_{i+2}^{k+1}x_{i+1}^{-1})^l
 x_{i+1}^{-k}x_{i+2}^{k+1})^{-1}$ and $i\in \mathbb{Z}/n$. Replacing $k$ and $l$ by $k-1$, $l-1$ respectively, we get 
 for each $i\in \mathbb{Z}/n$
 $$r_i=(x_i^{-k}x_{i+1}^{k})^l(x_{i+2}^{-k}x_{i+1}^{k})^{l-1}(
 x_{i+2}^{-k}x_{i+1}^{k-1})$$
Since there is an automorphism of $G$ given by sending $x_i$ to $x_{i+1}$, for each $i\in \mathbb{Z}/n$
and $K_{[2k,-2l]}$ is nontrivial no $x_i$ is 
trivial.
Set $x=x_1$, $y=x_2$, $z=x_3$, $t=x_4$, $w=x_5$. We have
\begin{itemize}
 \item $r_0=xyztw=1$,
 \item $r_1=(x^{-k}y^{k})^l(z^{-k}y^{k})^{l-1}(
 z^{-k}y^{k-1})=1$
 \item $r_2=(y^{-k}z^{k})^l(t^{-k}z^{k})^{l-1}(
 t^{-k}z^{k-1})=1$
 \item $r_3=(z^{-k}t^{k})^l(w^{-k}t^{k})^{l-1}(
 w^{-k}t^{k-1})=1$
 \item $r_4=(t^{-k}w^{k})^l(x^{-k}w^{k})^{l-1}(
 x^{-k}w^{k-1})=1$
 \item $r_5=(w^{-k}x^{k})^l(y^{-k}x^{k})^{l-1}(
 y^{-k}x^{k-1})=1.$
\end{itemize}

We prove Theorem \ref{thm: main result3} by contradiction. Let $<$ be a left-order on $G$. Without loss of generality $x>1$. Here is the list
of possible signs for $y$, $z$, $t$, $w$.
\begin{enumerate}
 \item $x>1$, $y>1$, $z>1$, $t>1$, $w>1$;
 \item$x>1$, $y<1$, $z>1$, $t>1$, $w>1$;
 \item $x>1$, $y<1$, $z<1$, $t>1$, $w>1$;
 \item $x>1$, $y<1$, $z<1$, $t<1$, $w>1$;
 \item $x>1$, $y<1$, $z<1$, $t<1$, $w<1$;
 \item $x>1$, $y>1$, $z<1$, $t>1$, $w>1$;
 \item $x>1$, $y>1$, $z<1$, $t<1$, $w>1$;
 \item $x>1$, $y>1$, $z<1$, $t<1$, $w<1$;
 \item $x>1$, $y>1$, $z>1$, $t<1$, $w>1$;
 \item $x>1$, $y>1$, $z>1$, $t<1$, $w<1$;
 \item $x>1$, $y>1$, $z>1$, $t>1$, $w<1$;
 \item $x>1$, $y<1$, $z>1$, $t>1$, $w<1$;
 \item $x>1$, $y<1$, $z>1$, $t<1$, $w<1$;
 \item $x>1$, $y<1$, $z<1$, $t>1$, $w<1$;
 \item $x>1$, $y<1$, $z>1$, $t<1$, $w>1$;
 \item $x>1$, $y>1$, $z<1$, $t>1$, $w<1$.
\end{enumerate}
We can immediately rule out eleven of these possibilities,

\begin{table}[H]\label{Table: 1}
\centering
\begin{tabular}{ | m{2cm} | m{4cm} | } 
\hline
Possibility & Ruled out by  \\ 
\hline
1 & $r_0$ \\ 
\hline
2 & $r_1$  \\ 
\hline
5 & $r_5$  \\ 
\hline
6 & $r_2$  \\ 
\hline
9 & $r_3$   \\ 
\hline
11 & $r_4$   \\ 
\hline
12 & $r_1$   \\ 
\hline
13 & $r_1$   \\ 
\hline
14 & $r_3$   \\ 
\hline
15 & $r_1$   \\ 
\hline
16 & $r_2$   \\ 
\hline

\end{tabular}
\caption{}
\end{table}
The cases which are left to check are (3), (4), (7), (8) and (10).

The automorphism $\phi$ of $G$ which sends
$(x,y,z,t,w)$ to $(w,x,y,z,t)$ acts on $LO(G)$. Using $\phi$ the reader will verify that up to replacing an order by 
its opposite, an order of the form  (3), (4), (7), (8) and (10) exists if and only if one of type 
$x>1$, $y>1$, $z<1$, $t<1$, $w<1$ exists. Assume that $<$ satisfies these inequalities.

\begin{lemma}\label{lemma: 2.2}
 $x^{-k}y^{k}<1$, $x^{k}y^{-k}<1$.
\end{lemma}

\begin{proof}
 Since $z<1$ then $z^{-1}>1$. By relation $r_1$, since $y>1$, we have that $x^{-k}y^{k}<1$. \\
 By relation $r_5=(w^{-k}x^{k})^l(y^{-k}x^{k})^{l-1}(
 y^{-k}x^{k-1})=1$ we have,
 $$(w^{-k}x^{k}w^{-k}x^{k}...w^{-k}x^{k})(y^{-k}x^{k}y^{-k}x^{k}...y^{-k}x^{k})(y^{-k}x^{k-1})=1,$$ this implies
 $$w^{-k}(x^{k}w^{-k})^{l-1}(x^{k}y^{-k})^{l}
 x^{k-1}=1.$$ Since $w^{-1}>1$, and $x>1$ then $x^{k}y^{-k}<1$.

\end{proof}

\begin{lemma}\label{lemma: 2.3}
 $wx>1$ and $yz>1$.
\end{lemma}

\begin{proof}
 Assume by contradiction that $wx<1$, then $x<w^{-1}\leq w^{-k}$, which implies that $x^{-1}w^{-k}>1$.\\
 By $r_5=(w^{-k}x^{k})^l(y^{-k}x^{k})^{l-1}(
 y^{-k}x^{k-1})=1$, we have,
 $$(w^{-k}x^{k})(w^{-k}x^{k})^{l-1}(y^{-k}x^{k})^{l-1}(y^{-k}x^{k})
 x^{-1}=1,$$ which implies 
 $$x^{k}(w^{-k}x^{k})^{l-1}(y^{-k}x^{k})^{l}
 x^{-1}w^{-k}=1.$$ Since $w^{-1}>1$, and $x>1$ then $y^{-k}x^{k}<1$, which implies $x^{-k}y^{k}>1$ 
 this is a contradiction to Lemma \ref{lemma: 2.2}.
 
 Let show now that $yz>1$. By contradiction assume $yz<1$ then $z^{-1}y^{-1}>1$. By relation
 $$r_1=(x^{-k}y^{k})^l(z^{-k}y^{k})^{l-1}(
 z^{-k}y^{k-1})=1$$ we obtain $$(x^{-k}y^{k}x^{-k}y^{k}...x^{-k}y^{k})(z^{-k}y^{k})^{l-1}(
 z^{-k}y^{k-1})=1$$ which implies

 $$(y^{k}x^{-k})^{l-1}y^{k}(z^{-k}y^{k})^{l-1}
 z^{-k}y^{k-1}x^{-k}=1$$
 By Lemma \ref{lemma: 2.2} we obtain $x^{k}y^{-k}<1$, which implies $y^{-k}<x^{-k}$, by multiplying in both side by 
 $(y^{k}x^{-k})^{l-1}y^{k}(z^{-k}y^{k})^{l-1}
 z^{-k}y^{k-1}$ we have that,
 $$(y^{k}x^{-k})^{l-1}y^{k}(z^{-k}y^{k})^{l-1}
 z^{-k}y^{k-1}y^{-k}<1$$
 which implies
 $$(y^{k}x^{-k})^{l-1}y^{k}(z^{-k}y^{k})^{l-1}
 z^{-(k-1)}z^{-1}y^{-1}<1$$
 which is a contradiction since the terms on the left-hand side are positive. Thus $yz>1$.

\end{proof}

\begin{lemma}\label{lemma: 2.4}
 $z^{k}t^{-k}>1$, $t^{-k}z^{k}>1$, and $t^{k-1}>w^{k}>t^{k}$.
\end{lemma}

\begin{proof}

By $r_2=(y^{-k}z^{k})^l(t^{-k}z^{k})^{l-1}(
 t^{-k}z^{k-1})=1$ we have $$(y^{-k}z^{k}y^{-k}z^{k}...y^{-k}z^{k})(t^{-k}z^{k})^{l-1}(
 t^{-k}z^{k-1})=1$$ which implies $$ y^{-k}(z^{k}y^{-k})^{l-1}(z^{k}t^{-k})^{l}
 z^{k-1}=1.$$ Since $y^{-1}<1$ and $z<1$, then $z^{k}t^{-k}>1$.\\
Again by  $r_2=(y^{-k}z^{k})^l(t^{-k}z^{k})^{l-1}(
 t^{-k}z^{k-1})=1$  we obtain, $$y^{-k}z^{k}(y^{-k}z^{k})^{l-1}(t^{-k}z^{k})^{l}z^{-1}=1$$
which implies $$z^{k}(y^{-k}z^{k})^{l-1}(t^{-k}z^{k})^{l}z^{-1}y^{-k}=1.$$
Since $yz>1$ by Lemma \ref{lemma: 2.3}, then $z>y^{-1}\geq y^{-k}$, which implies $1>z^{-1}y^{-k}$.
Therefore, $t^{-k}z^{k}>1$.

Now it left to show that $t^{k-1}>w^{k}>t^{k}$. By $$r_4=(t^{-k}w^{k})^l(x^{-k}w^{k})^{l-1}(
 x^{-k}w^{k-1})=1$$ since $x^{-1}<1$ and $w<1$, then $t^{-k}w^{k}>1$, which implies that $t^{k}<w^{k}$ and $w^{-k}t^{k}<1$.
Since $t^{-k}z^{k}>1$ then $z^{-k}t^{k}<1$, and by relation 
$$r_3=(z^{-k}t^{k})^l(w^{-k}t^{k})^{l-1}(
 w^{-k}t^{k-1})=1$$  we obtain $w^{-k}t^{k-1}>1$. Thus, $t^{k-1}>w^{k}$.

\end{proof}
\begin{lemma}\label{lemma: 2.5}
$t^{-1}>w^{-k+1}x^k$ and $w^{-k}<t^{-k}$.
\end{lemma}
\begin{proof}
The equation  $$r_3=(z^{-k}t^{k})^l(w^{-k}t^{k})^{l-1}(
 w^{-k}t^{k-1})=1$$ implies  $$A=t^{-1}(z^{-k}t^{k})^l(w^{-k}t^{k})^{l}=1$$ multiplying $A$ with 
 $$r_4=(t^{-k}w^{k})^l(x^{-k}w^{k})^{l-1}(x^{-k}w^{k-1})=1$$ we obtain  
 $$Ar_4=t^{-1}(z^{-k}t^{k})^l(x^{-k}w^{k})^{l-1}(x^{-k}w^{k-1})=1$$
 which implies
 $$(z^{-k}t^{k})^l(x^{-k}w^{k})^{l-1}(x^{-k}w^{k-1}t^{-1})=1.$$
 Since $z^{-k}t^{k}<1$ by Lemma \ref{lemma: 2.4}, and $x^{-k}w^{k}<1$ because $x^{-1}<1$ and $w<1$, then $x^{-k}w^{k-1}t^{-1}>1$.
 Therefore $t^{-1}>w^{-k+1}x^k$.
 
Now, let show that $w^{-k}<t^{-k}$. The equation $$r_5=(w^{-k}x^{k})^l(y^{-k}x^{k})^{l-1}(y^{-k}x^{k-1})=1$$ implies 
$$w^{-k}x^{k}(w^{-k}x^{k})^{l-1}(y^{-k}x^{k})^{l}x^{-1}=1$$ which implies 
$$x^{-1}w^{-k}x^{k}(w^{-k}x^{k})^{l-1}(y^{-k}x^{k})^{l}=1$$ Now, $w^{-k}x^{k}>1$, and by Lemma \ref{lemma: 2.2}
$y^{-k}x^{k}>1$. Therefore
$x^{-1}w^{-k}x^{k}<1$. This implies $1<x^{k}<w^{k}x$. Hence $x>w^{-k}$ , $x^{k}>x>w^{-k}>w^{-1}$ and  $x^{k}>w^{-1}$.
Multiplying both sides of $x^{k}>w^{-1}$ by $w^{-k+1}$ we have $w^{-k+1}x^{k}>w^{-k+1}w^{-1}=w^{-k}$, and since
$t^{-1}>w^{-k+1}x^k$ then
$t^{-k}>t^{-1}>w^{-k+1}x^k>w^{-k}$. Therefore $w^{-k}<t^{-k}$.

\end{proof}
\begin{remark}
{\rm $\;$ 
 Lemmas \ref{lemma: 2.2}, \ref{lemma: 2.3}, \ref{lemma: 2.4}, \ref{lemma: 2.5} hold for all lo's for
 which $1 < x, 1 < y, z < 1, t < 1, w < 1$. In particular, they will apply to $<_c$ in the Proof of 
 Teorem \ref{thm: main result3}. 
 }
\end{remark}

Since $G$ is a countable left-orderable group, it can be seen as a subgroup of $Homeo_+(\mathbb{R})$ and it
will act effectively on $\mathbb{R}$ by order preserving homeomorphisms without global fixed points.\\

\begin{lemma} \label{lemma: 2.6}
Let $a$ be a fixed point for $x$. Then the only possible signs for $x$, $y$, $z$, $t$, $w$ with respect to $<_a$ are  
$t<_a1$, $w<_a1$, $1<_ax$, $1<_ay$, $1<_az$.
\end{lemma}
\begin{proof}
 Since $x>1$ and $x(a)=a$, we have $1<_ax$. Therefore as in the analysis of (1) through (16) above, 
 one can see that the cases which are left to check are
 \begin{enumerate}
 \item $y<_a1$, $z<_a1$, $1<_at$, $1<_aw$, $1<_ax$
 \item $1<_aw$, $1<_ax$, $y<_a1$, $z<_a1$, $t<_a1$
 \item $z<_a1$, $t<_a1$, $1<_aw$, $1<_ax$, $1<_ay$
 \item $1<_ax$, $1<_ay$, $z<_a1$, $t<_a1$, $w<_a1$
 \item $t<_a1$, $w<_a1$, $1<_ax$, $1<_ay$, $1<_az$
\end{enumerate}
For (1), let $<'=<_a^{op}$, using Lemma \ref{lemma: 2.4} and replacing $x$ by $y$, $y$ by $z$, $z$ by $t$, $t$ by $w$ and
$w$ by $x$, we have $w^{k}<'x^{k}<'w^{k-1}$. Therefore  $w^{k-1}<_ax^{k}<_aw^{k}$, which
implies $w(a)=a$, which is impossible because $w^{-1}(a)=a$ and $w^{-1}>1$ implies that $1<_aw^{-1}$, which is a contradiction 
to the fact that $1<_aw$. Therefore case (1) is not possible. 

For (2), using Lemma \ref{lemma: 2.2} and replacing $x$ by $w$, $y$ by $x$, $z$ by $y$, $t$ by $z$ and
$w$ by $t$, we have $x^{k}<_aw^{k}$, $x^{-k}<_aw^{-k}$. Hence $w(a)=a$ and by the same argument as for (1) we get a 
contradiction.  Therefore case (2) is not possible.

For (3), let $<'=<_a^{op}$, using Lemma \ref{lemma: 2.4} and replacing $x$ by $z$, $y$ by $t$, $z$ by $w$, $t$ by $x$ and
$w$ by $y$, we have $x^{k}<'w^{k}$, hence $w^k<_ax^k$ and $w^k(a)\leq a$, and since $1<_aw$ then $w(a)=a$
and by the same argument as for (1) we get a 
contradiction.  Therefore case (3) is not possible.

For (4), Lemma \ref{lemma: 2.2} implies $x^{-k}y^{k}(a)\leq a$, $x^{k}y^{-k}(a)\leq a$, hence 
$a=x^{-k}(a)\geq y^{-k}(a)$ and  $y^{k}(a)\leq x^{k}(a)=a$. This two inequalities imply that $y(a)=a$. Therefore, by relation
$$r_5=(w^{-k}x^{k})^l(y^{-k}x^{k})^{l-1}(
 y^{-k}x^{k-1})=1$$ we obtain that $(w^{-k}x^{k})^l(a)=a$, which implies $(w^{-k}x^{k})(a)=a$. Hence $w^{-k}(a)=a$, and 
 $w(a)=a$. A similar argument using relation $r_4$ shows that $t(a)=a$. Since $xyztw=1$, $z(a)=a$.
 Thus $a$ is fixed by $G$, which contradicts our assumptions.  Therefore case (4) is not possible.
\end{proof}

\begin{lemma}\label{lemma: 2.7}
Let $b$ be a fixed point for $y$. Then the only possible signs for $x$, $y$, $z$, $t$, $w$ with respect to $<_b$ are  
$w<_b1$, $x<_b1$, $1<_by$, $1<_bz$, $1<_bt$.
\end{lemma}
\begin{proof}
Since $y>1$ and $y(b)=b$, we have $1<_by$. Therefore as in the analysis of (1) through (16) above, by replacing
$x$ by $y$, $y$ by $z$, $z$ by $t$, $t$ by $w$, $w$ by $x$ and $<$ by $<_b$, one can see that the cases which are
left to check are 
 \begin{enumerate}
 \item $t<_b1$, $w<_b1$, $1<_bx$, $1<_by$, $1<_bz$
 \item $1<_by$, $1<_bz$, $t<_b1$, $w<_b1$, $x<_b1$
 \item $z<_b1$, $t<_b1$, $1<_bw$, $1<_bx$, $1<_by$
 \item $1<_bx$, $1<_by$, $z<_b1$, $t<_b1$, $w<_b1$
 \item $w<_b1$, $x<_b1$, $1<_by$, $1<_bz$, $1<_bt$.
\end{enumerate}
For (1), let $<'=<_b^{op}$, using Lemma \ref{lemma: 2.4} and replacing $x$ by $t$, $y$ by $w$, $z$ by $x$, $t$ by $y$ and
$w$ by $z$, we have  $y^{k}<'z^{k}<'y^{k-1}$. Therefore  $y^{k-1}<_bz^{k}<_bw^{k}$, which
implies $z(b)=b$, which is impossible because $z^{-1}(b)=b$ and $z^{-1}>1$ implies that $1<_bz^{-1}$, which is a contradiction 
to the fact that $1<_bz$.  Therefore case (1) is not possible.

For (2), using Lemma \ref{lemma: 2.2} and replacing $x$ by $y$, $y$ by $z$, $z$ by $t$, $t$ by $w$ and
$w$ by $x$, we have $z^{k}<_by^{k}$, $z^{-k}<_by^{-k}$. Hence $z(b)=b$ and by the same argument as for (1) we get a 
contradiction. Therefore case (2) is not possible.

For (3), let $<'=<_b^{op}$, using Lemma \ref{lemma: 2.4} and replacing $x$ by $z$, $y$ by $t$, $z$ by $w$, $t$ by $x$ and
$w$ by $y$, we have $x^{k}<'y^{k}<'x^{k-1}$. Therefore  $x^{k-1}<_by^{k}<_bx^{k}$, which
implies $x(b)=b$. By equation $$r_5=(w^{-k}x^{k})^l(y^{-k}x^{k})^{l-1}(
 y^{-k}x^{k-1})=1$$ we obtain that $(w^{-k}x^{k})^l(b)=b$, which implies $(w^{-k}x^{k})(b)=b$. Hence $w^{-k}(b)=b$, and 
 $w(b)=b$. This is impossible because $w^{-1}(b)=b$ and $w^{-1}>1$ implies that $1<_bw^{-1}$, which is a contradiction 
to the fact that $1<_bw$. Therefore case (3) is not possible.

For (4), Lemma \ref{lemma: 2.2} implies $x^{-k}y^{k}(b)\leq b$, $x^{k}y^{-k}(b)\leq b$, hence 
$x^{-k}(b)\geq y^{-k}(b)=b$ and  $b=y^{k}(b)\leq x^{k}(b)$. This two inequalities imply that $x(b)=b$. Therefore, by relation
$$r_5=(w^{-k}x^{k})^l(y^{-k}x^{k})^{l-1}(
 y^{-k}x^{k-1})=1$$ we obtain that $(w^{-k}x^{k})^l(b)=b$, which implies $(w^{-k}x^{k})(b)=b$. Hence $w^{-k}(b)=b$, and 
 $w(b)=b$. A similar argument using relation $r_4$ shows that $t(b)=b$. Since $xyztw=1$, $z(b)=b$.
 Thus $b$ is fixed by $G$, which contradicts our assumptions. Therefore case (4) is not possible.
\end{proof}

\begin{lemma}\label{lemma: 2.8}

Let $c$ be a fixed point for $z$. Then the only possible signs for $x$, $y$, $z$, $t$, $w$ with respect to $<_c$ are  
$1<_cx$, $1<_cy$, $z<_c1$, $t<_c1$, $w<_c1$.
\end{lemma}
\begin{proof}
Since $z<1$ and $z(c)=c$, we have $z<_c1$. Let  $<'=<_c^{op}$, then $1<'z$.
Therefore as in the analysis of (1) through (16) above, by replacing
$x$ by $z$, $y$ by $t$, $z$ by $w$, $t$ by $x$, $w$ by $y$ and $<$ by $<'$,
 one can see that the cases which are left to check are 
 \begin{enumerate}
 \item $1<'z$, $1<'t$, $w<'1$, $x<'1$, $y<'1$
 \item $t<'1$, $w<'1$, $1<'x$, $1<'y$, $1<'z$
 \item $w<'1$, $x<'1$, $1<'y$, $1<'z$, $1<'t$
 \item $x<'1$, $y<'1$, $1<'z$, $1<'t$, $1<'w$
 \item $1<'y$, $1<'z$, $t<'1$, $w<'1$, $x<'1$
\end{enumerate}
Therefore,
 \begin{enumerate}
 \item $z<_c1$, $t<_c1$, $1<_cw$, $1<_cx$, $1<_cy$
 \item $1<_ct$, $1<_cw$, $x<_c1$, $y<_c1$, $z<_c1$
 \item $1<_cw$, $1<_cx$, $y<_c1$, $z<_c1$, $t<_c1$
 \item $1<_cx$, $1<_cy$, $z<_c1$, $t<_c1$, $w<_c1$
 \item $y<_c1$, $z<_c1$, $1<_ct$, $1<_cw$, $1<_cx$
\end{enumerate}
For (1), let $<'=<_c^{op}$, using Lemma \ref{lemma: 2.2} and replacing $x$ by $z$, $y$ by $t$, $z$ by $w$, $t$ by $x$ and
$w$ by $y$, we have $t^{k}<'z^{k}$, $t^{-k}<'z^{-k}$. Therefore $z^{k}<_ct^{k}$, $z^{-k}<_ct^{-k}$, which 
implies $t(c)=c$, using equation $r_2$ we have that $y(c)=c$ and using equation $r_1$ that $x(c)=c$, and therefore by equation
$r_5$, we have that $w(c)=c$, which is impossible because $w^{-1}(c)=c$ and $w^{-1}>1$ implies that $1<_cw^{-1}$,
which is a contradiction 
to the fact that $1<_cw$. Therefore case (1) is not possible. 

For (2), using Lemma \ref{lemma: 2.4} and replacing $x$ by $t$, $y$ by $w$, $z$ by $x$, $t$ by $y$ and
$w$ by $z$ we have $y^{k}<_cz^{k}<_cy^{k-1}$, which
implies $y(c)=c$. This is impossible because $y(c)=c$ and $y>1$ implies that $1<_cy$, which is a contradiction 
to the fact that $y<_c1$. Therefore case (2) is not possible.

For (3), using Lemma \ref{lemma: 2.4} and replacing $x$ by $w$, $y$ by $x$, $z$ by $y$, $t$ by $z$ and
$w$ by $t$, we have $z^{k}<_ct^{k}<_cz^{k-1}$, which
implies $t(c)=c$, using equation $r_2$ we have that $y(c)=c$,  which is impossible because $y(c)=c$ and $y>1$ 
implies that $1<_cy$, which is a contradiction 
to the fact that $y<_c1$.
Therefore case (3) is not possible.

For (5), let $<'=<_c^{op}$, using Lemma \ref{lemma: 2.2}, and replacing $x$ by $y$, $y$ by $z$, $z$ by $t$, $t$ by $w$ and
$w$ by $x$, we have $z^{k}<'y^{k}$, $z^{-k}<'y^{-k}$. Therefore $y^{k}<_cz^{k}$, $y^{-k}<_cz^{-k}$, which 
implies $y(c)=c$. This is impossible because $y(c)=c$ and $y>1$ 
implies that $1<_cy$, which contradicts 
the fact that $y<_c1$. Therefore case (5) is not possible.
\end{proof}

\begin{lemma}\label{lemma: 2.9}
Let $d$ be a fixed point for $t$. Then the only possible signs for $x$, $y$, $z$, $t$, $w$ with respect to $<_d$ are  
$1<_dy$, $1<_dz$, $t<_d1$, $w<_d1$, $x<_d1$.
\end{lemma}
\begin{proof}
 Similar proof as for Lemma \ref{lemma: 2.6}.

\end{proof}

\begin{lemma}\label{lemma: 2.10}
Let $e$ be a fixed point for $w$. Then the only possible signs for $x$, $y$, $z$, $t$, $w$ with respect to $<_e$ are  
$1<_ez$, $1<_et$, $w<_e1$, $x<_e1$, $y<_e1$.
\end{lemma}

\begin{proof}
 Similar proof as for Lemma \ref{lemma: 2.6}.
\end{proof}
Summarizing Lemmas \ref{lemma: 2.6} through \ref{lemma: 2.10}, we have
\begin{remark}
{\rm $\;$ 

(1) If $e$ is a fixed point for $w$, then the only possible signs for $x$, $y$, $z$, $t$, $w$ with respect to $<_e$ are  
$1<_ez$, $1<_et$, $w<_e1$, $x<_e1$, $y<_e1$. Therefore $x$, $y$, $z$, $t$ will also have fixed points (different to $e$).

(2) If $d$ is a fixed point for $t$, then the only possible signs for $x$, $y$, $z$, $t$, $w$ with respect to $<_d$ are  
$1<_dy$, $1<_dz$, $t<_d1$, $w<_d1$, $x<_d1$. Therefore $x$, $z$ will also have fixed points (different to $d$).

(3) If  $b$ is a fixed point for $y$, then the only possible signs for $x$, $y$, $z$, $t$, $w$ with respect to $<_b$ are  
$w<_b1$, $x<_b1$, $1<_by$, $1<_bz$, $1<_bt$. Therefore $x$, $z$, $t$ will also have fixed points (different to $b$).

(4) If $a$ is a fixed point for $x$, then the only possible signs for $x$, $y$, $z$, $t$, $w$ with respect to $<_a$ are  
$t<_a1$, $w<_a1$, $1<_ax$, $1<_ay$, $1<_az$. Therefore $z$ also has a fixed point (different to $a$).
}
\end{remark}

\begin{lemma}\label{lemma: 2.11}
 $z$ has a fixed point.
\end{lemma}

\begin{proof}
 Assume by contradiction that $z$ is fixed point free. Then by the previous remark, $x$, $y$, $t$, $w$ are also fixed point free.
 Therefore, for any $a\in \mathbb{R}$ the only possible signs for $x$, $y$, $z$, $t$, $w$ with respect to $<_a$ are  
$1<_ax$, $1<_ay$, $z<_a1$, $t<_a1$, $w<_a1$. Then by Lemma \ref{lemma: 2.2},
$x^{-k}y^{k}(a)\leq a$, $x^{k}y^{-k}(a)\leq a$, which implies
$y^{k}(a)\leq x^{k}(a)$ and $y^{-k}(a)\leq x^{-k}(a)$. Since $a\in \mathbb{R}$ was arbitrary, the last two 
inequalities are true for any $a\in \mathbb{R}$. Then $a\leq y^{-k}(x^{k}(a))$ for any $a\in \mathbb{R}$, which implies
$a\leq y^{-k}(x^{k}(a))\leq  x^{-k}(x^{k}(a))=a$, hence $x^{k}(a)=y^{k}(a)$ for any $a\in \mathbb{R}$. Therefore
$x^{k}=y^{k}$. We obtain a contradiction from the equation
$$r_1=(x^{-k}y^{k})^l(z^{-k}y^{k})^{l-1}(
 z^{-k}y^{k-1})=1$$ since $y>1$, $z^{-1}>1$ and $x^{-k}y^{k}=1$.
\end{proof}

\begin{proof}[Proof of Theorem \ref {thm: main result3}]

The equation  $$r_3=(z^{-k}t^{k})^l(w^{-k}t^{k})^{l-1}(
 w^{-k}t^{k-1})=1$$ implies $$z^{-k}(t^{k}z^{-k})^{l-1}(t^{k}w^{-k})^{l}t^{k-1}=1$$ which implies 
 $$t^{k-1}z^{-k}(t^{k}z^{-k})^{l-1}(t^{k}w^{-k})^{l}=1.$$ Since $t^{k}z^{-k}<1$ by Lemma \ref{lemma: 2.4},
 $t^{k}w^{-k}<1$ by Lemma \ref{lemma: 2.5}, and $t^{k-1}z^{-k}(t^{k}z^{-k})^{l-1}(t^{k}w^{-k})^{l}=1$, 
 then $t^{k-1}z^{-k}\geq 1$.
 
By Lemma \ref{lemma: 2.11} $z$ has a fixed point, call it $c$.
By Lemma \ref{lemma: 2.8}, the only possible signs for $x$, $y$, $z$, $t$, $w$ with respect to $<_c$ are  
$1<_cx$, $1<_cy$, $z<_c1$, $t<_c1$, $w<_c1$. Then $1\leq_c t^{k-1}z^{-k}$.
This implies that  $t^{-k+1}(c)\leq z^{-k}(c)=c$. 
Hence $t^{-k+1}(c)=c$, which implies $t(c)=c$, and this is impossible by Lemma \ref{lemma: 2.9}.
\end{proof}
\section{Left-orderability and genus 2 two-bridge knots}
In this section we will prove the following theorem.
\begin{theorem}\label{thm: main result2}
  The fundamental group of the 3-fold cyclic branched cover of $K_{[-2q,2s,-2t,2l]}$ is not left-orderable
  where $q$, $s$, $t$ and $l$ 
  $\in\mathbb{Z}\setminus\{0\}$.
\end{theorem}

Mulazzani and Vesnin ([MV], Theorem 8) proved that the generalized periodic Takahashi manifold 
$T_{n,m}(\frac{1}{a_j}; \frac{1}{b_j})$ is 
the $n$-fold cyclic branched covering of the two-bridge knot corresponding to the Conway parameters 
$[-2a_1, 2b_1, -2a_2,2 b_2,  \cdots, -2a_m, 2b_m]$.
For $m=2$ we have the family of genus 2 two-bridge knots $K_{[-2q,2s,-2t,2l]}$. In order to show that the 
fundamental group of the 3-fold cyclic branched cover of any genus 2 two-bridge knot $K_{[-2q,2s,-2t,2l]}$ 
is not left-orderable we have to consider all cases for the signs of $q$, $s$, $t$ and $l$. We have sixteen
cases to consider, but since the mirror image of the knot $K_{[-2q,2s,-2t,2l]}$ is the knot $K_{[2q,-2s,2t,-2l]}$ we need only
deal with eight of them:
\begin{enumerate}
 \item $q>0$, $s>0$, $t>0$ and $l>0$;
 \item $q>0$, $s>0$, $t<0$ and $l>0$;
 \item $q<0$, $s>0$, $t<0$ and $l>0$;
 \item $q<0$, $s<0$, $t<0$ and $l>0$;
 \item $q<0$, $s>0$, $t>0$ and $l>0$;
 \item $q<0$, $s>0$, $t<0$ and $l<0$;
 \item $q>0$, $s<0$, $t<0$ and $l>0$;
 \item $q>0$, $s>0$, $t<0$ and $l<0$.
 
\end{enumerate}

By Mulazzani and Vesnin ([MV], Theorem 10), the fundamental group of the $n$-fold cyclic branched cover of $K_{[-2q,2s,-2t,2l]}$ is 
$$G:=\pi_1(\Sigma_n(K_{[-2q,2s,-2t,2l]}))=\langle x_1, x_2,\cdots, x_n \;|\; w(x_{k-2}, x_{k-1},x_{k},x_{k+1},x_{k+2}),
\; k\in \mathbb{Z}/n\rangle$$
where\\
$w(x_{k-2}, x_{k-1},x_{k},x_{k+1},x_{k+2})=$

$[[(x_k^qx_{k+1}^{-q})^{-s}x_k(x_{k-1}^qx_{k}^{-q})^{s}]^tx_k^qx_{k+1}^{-q} 
[(x_{k+1}^qx_{k+2}^{-q})^{-s}x_{k+1}(x_{k}^qx_{k+1}^{-q})^{s}]^{-t}]^{-l}$
$(x_k^qx_{k+1}^{-q})^{-s}x_k(x_{k-1}^qx_{k}^{-q})^{s}$\\
$\cdot[[(x_{k-1}^qx_{k}^{-q})^{-s}x_{k-1}(x_{k-2}^qx_{k-1}^{-q})^{s}]^tx_{k-1}^qx_{k}^{-q} 
[(x_{k}^qx_{k+1}^{-q})^{-s}x_{k}(x_{k-1}^qx_{k}^{-q})^{s}]^{-t}]^{l}$ and $k\in \mathbb{Z}/n$.

Since there is an automorphism of $G$ given by sending $x_k$ to $x_{k+1}$, for each $k\in \mathbb{Z}/n$ and 
$K_{[-2q,2s,-2t,2l]}$ is nontrivial no $x_k$ is trivial.

For $n=3$, let $x=x_1$, $y=x_2$ and $z=x_3$, we have

 $r_1=[[(y^qx^{-q})^{s}x(z^qx^{-q})^{s}]^tx^qy^{-q} 
[(z^qy^{-q})^{s}y(x^qy^{-q})^{s}]^{-t}]^{-l}(y^qx^{-q})^{s}x(z^qx^{-q})^{s}$\\
$\cdot[[(x^qz^{-q})^{s}z(y^qz^{-q})^{s}]^tz^qx^{-q} 
[(y^qx^{-q})^{s}x(z^qx^{-q})^{s}]^{-t}]^{l}=1$;

 $r_2=[[(z^qy^{-q})^{s}y(x^qy^{-q})^{s}]^ty^qz^{-q} 
[(x^qz^{-q})^{s}z(y^qz^{-q})^{s}]^{-t}]^{-l}(z^qy^{-q})^{s}y(x^qy^{-q})^{s}$\\
$\cdot[[(y^qx^{-q})^{s}x(z^qx^{-q})^{s}]^tx^qy^{-q} 
[(z^qy^{-q})^{s}y(x^qy^{-q})^{s}]^{-t}]^{l}=1$;

$r_3=[[(x^qz^{-q})^{s}z(y^qz^{-q})^{s}]^tz^qx^{-q} 
[(y^qx^{-q})^{s}x(z^qx^{-q})^{s}]^{-t}]^{-l}(x^qz^{-q})^{s}z(y^qz^{-q})^{s}$\\
$\cdot[[(z^qy^{-q})^{s}y(x^qy^{-q})^{s}]^ty^qz^{-q} 
[(x^qz^{-q})^{s}z(y^qz^{-q})^{s}]^{-t}]^{l}=1$.

\begin{proof}[Proof of Theorem \ref {thm: main result2}]
 
Considering the product $r_3r_2r_1$ we have, $$zyx=1$$
We prove Theorem \ref {thm: main result2} by contradiction. Let $<$ be a left-order on $G$. Without loss of generality $x>1$. Here is the list
of possible signs for $y$, $z$.
\begin{enumerate}
 \item $x>1$, $y>1$, $z>1$;
 \item$x>1$, $y<1$, $z>1$;
 \item $x>1$, $y<1$, $z<1$;
 \item $x>1$, $y>1$, $z<1$;
\end{enumerate}

We can immediately rule out $(1)$ because $zyx=1$.
The cases which are left to check are (2), (3), (4).

The automorphism $\phi$ of $G$ which sends
$(x,y,z)$ to $(z,x,y)$ acts on $LO(G)$. Using $\phi$ the reader will verify that up to replacing an order by 
its opposite, an order of the form  (2), (3), (4),  exists if and only if one of type 
$x>1$, $y<1$, $z<1$ exists. Assume that $<$ satisfies these inequalities.

Let $$X=(y^qx^{-q})^{s}x(z^qx^{-q})^{s}$$ $$Y=(z^qy^{-q})^{s}y(x^qy^{-q})^{s}$$ and
$$Z=(x^qz^{-q})^{s}z(y^qz^{-q})^{s}.$$
Considering the product $r_3r_2r_1$ we have,
$ZYX=1$ and $r_1$, $r_2$, and $r_3$ become,
\begin{itemize}
 \item $r'_1=(Y^ty^qx^{-q}X^{-t})^{l}X(Z^tz^qx^{-q}X^{-t})^{l}=1$;

\item $r'_2=(Z^tz^qy^{-q}Y^{-t})^{l}Y(X^tx^qy^{-q}Y^{-t})^{l}=1$;

\item $r'_3=(X^tx^qz^{-q}Z^{-t})^{l}Z(Y^ty^qz^{-q}Z^{-t})^{l}=1$.

\end{itemize}

This is the same as,
\begin{itemize}
 \item $r''_1=(Y^ty^qx^{-q}X^{-t})^{l-1}Y^ty^qx^{-q}X^{-t+1}(Z^tz^qx^{-q}X^{-t})^{l}=1$;

\item $r''_2=(Z^tz^qy^{-q}Y^{-t})^{l-1}Z^tz^qy^{-q}Y^{-t+1}(X^tx^qy^{-q}Y^{-t})^{l}=1$;

\item $r''_3=(X^tx^qz^{-q}Z^{-t})^{l-1}X^tx^qz^{-q}Z^{-t+1}(Y^ty^qz^{-q}Z^{-t})^{l}=1$.

\end{itemize}
We have eight cases:

\textbf{First case:} $q>0$, $s>0$, $t>0$, $l>0$.

We have $X=(y^qx^{-q})^{s}x(z^qx^{-q})^{s}=(y^qx^{-q})^{s-1}y^qx^{-q+1}(z^qx^{-q})^{s}$. Since 
$y<1$, $z<1$ and $x^{-1}<1$ then $X<1$. 

We discuss the signs of $Y$ and $Z$. We have three subcases
\begin{enumerate}
 \item If $Y<1$, then since $ZYX=1$ and $X<1$ then $Z>1$. Since 
 $Y=(z^qy^{-q})^{s}y(x^qy^{-q})^{s}=(z^qy^{-q})^{s}yxx^{q-1}y^{-q}(x^qy^{-q})^{s-1}
 =(z^qy^{-q})^{s}z^{-1}x^{q-1}y^{-q}(x^qy^{-q})^{s-1}$
then $z^qy^{-q}<1$. We have,
 \begin{equation} \label{eq1}
\begin{split}
 X^tx^qz^{-q} &= X^{t-1}(y^qx^{-q})^{s}x(z^qx^{-q})^{s}x^qz^{-q}\\
 &= X^{t-1}(y^qx^{-q})^{s}x(z^qx^{-q})^{s-1}z^qx^{-q}x^qz^{-q}\\
 &= X^{t-1}(y^qx^{-q})^{s}x(z^qx^{-q})^{s-1}\\
 &= X^{t-1}(y^qx^{-q})^{s-1}y^qx^{-q+1}(z^qx^{-q})^{s-1}
\end{split}
\end{equation}
Since $X<1$, $y<1$, $z<1$ and $x^{-1}<1$ then $X^tx^qz^{-q}<1$. Therefore
$X^tx^qz^{-q}Z^{-t}<1$ and $X^tx^qz^{-q}Z^{-t+1}<1$.

We have also that,
\begin{equation} \label{eq1}
\begin{split}
y^qz^{-q}Z^{-t} &= y^qz^{-q}(z^qy^{-q})^{s}z^{-1}(z^qx^{-q})^{s}Z^{-t+1}\\
 &= y^qz^{-q}z^qy^{-q}(z^qy^{-q})^{s-1}z^{-1}(z^qx^{-q})^{s}Z^{-t+1}\\
 &= (z^qy^{-q})^{s-1}z^{q-1}x^{-q}(z^qx^{-q})^{s-1}Z^{-t+1}
\end{split}
\end{equation}
Since $Z^{-1}<1$, $z^qy^{-q}<1$, $z<1$ and $x^{-1}<1$ then $y^qz^{-q}Z^{-t}<1$. Therefore
$Y^ty^qz^{-q}Z^{-t}<1$, and we have a contradiction by the relation $r''_3$.
 
 \item If $Z<1$ then since $ZYX=1$ and $X<1$ then $Y>1$. 
 Since $Z=(x^qz^{-q})^{s}z(y^qz^{-q})^{s}=(x^qz^{-q})^{s-1}x^qz^{-q+1}(y^qz^{-q})^{s}$ then
 $y^qz^{-q}<1$. We have
\begin{equation} \label{eq1}
\begin{split}
z^qy^{-q}Y^{-t} &= z^qy^{-q}(y^qx^{-q})^{s}y^{-1}(y^qz^{-q})^{s}Y^{-t+1}\\
 &= z^qy^{-q}y^qx^{-q}(y^qx^{-q})^{s-1}y^{-1}(y^qz^{-q})^{s}Y^{-t+1}\\
&= z^qx^{-q}(y^qx^{-q})^{s-1}y^{-1}(y^qz^{-q})^{s}Y^{-t+1}\\
 &= z^qx^{-q}(y^qx^{-q})^{s-2}y^{q}x^{-q+1}x^{-1}y^{-1}(y^qz^{-q})^{s}Y^{-t+1}\\
 &= z^qx^{-q}(y^qx^{-q})^{s-2}y^{q}x^{-q+1}z(y^qz^{-q})^{s}Y^{-t+1}
\end{split}
\end{equation}
Since
$Y^{-1}<1$ and $y^qz^{-q}<1$ then $z^qy^{-q}Y^{-t}<1$. Therefore $Z^tz^qy^{-q}Y^{-t}<1$ and $Z^tz^qy^{-q}Y^{-t+1}<1$ for $t>1$.
For $t=1$, the relation $r''_2$ become 
$$r'''_2=(Zz^qy^{-q}Y^{-1})^{l-1}Zz^qy^{-q}Xx^qy^{-q}Y^{-1}(Xx^qy^{-q}Y^{-1})^{l-1}=1$$ and we have
$$z^qy^{-q}X=z^qy^{-q}y^qx^{-q}(y^qx^{-q})^{s-1}x(z^qx^{-q})^{s}=z^qx^{-q}(y^qx^{-q})^{s-2}y^qx^{-q+1}(z^qx^{-q})^{s}<1.$$

We have also that,
\begin{equation} \label{eq1}
\begin{split}
Y^ty^qx^{-q} &= Y^{t-1}(z^qy^{-q})^{s}y(x^qy^{-q})^{s}y^qx^{-q}\\
 &= Y^{t-1}(z^qy^{-q})^{s}y(x^qy^{-q})^{s-1}x^qy^{-q}y^qx^{-q}\\
 &= Y^{t-1}(z^qy^{-q})^{s}y(x^qy^{-q})^{s-1}= Y^{t-1}(z^qy^{-q})^{s}yxx^{q-1}y^{-q}(x^qy^{-q})^{s-2}\\
 &= Y^{t-1}(z^qy^{-q})^{s}z^{-1}x^{q-1}y^{-q}(x^qy^{-q})^{s-2}
\end{split}
\end{equation}
Since $z^qy^{-q}>1$, $z^{-1}>1$ and $y^{-1}>1$ then $Y^ty^qx^{-q}>1$. Hence 
$Y^ty^qx^{-q}X^{-t}>1$ which implies $X^tx^qy^{-q}Y^{-t}=(Y^ty^qx^{-q}X^{-t})^{-1}<1$. Therefore by relation
$r''_2$ and $r'''_2$ 
we have
a contradiction.
\item The only subcase which is left to check is if $Y>1$ and $Z>1$. In this subcase we have two subsubcases
\begin{itemize}
 \item If $z^qy^{-q}>1$ then by the last part of subcase (2) we have $Y^ty^qx^{-q}X^{-t}>1$, $Y^ty^qx^{-q}X^{-t+1}>1$ and
by the first part of subcase (1) we have that 
 $X^tx^qz^{-q}Z^{-t}<1$ which implies $Z^tz^qx^{-q}X^{-t}=(X^tx^qz^{-q}Z^{-t})^{-1}>1$. This gives a contradiction by the 
 relation $r''_1$.
 \item If $z^qy^{-q}<1$ then $y^qz^{-q}>1$. We have
 \begin{equation} \label{eq1}
\begin{split}
y^qx^{-q}X^{-t} &= y^qx^{-q}(x^qz^{-q})^{s}x^{-1}(x^qy^{-q})^{s}X^{-t+1}\\
 &= y^qx^{-q}x^qz^{-q}(x^qz^{-q})^{s-1}x^{-1}(x^qy^{-q})^{s}X^{-t+1}\\
 &= y^qz^{-q}(x^qz^{-q})^{s-1}x^{-1}(x^qy^{-q})^{s}X^{-t+1}\\
&= y^qz^{-q}(x^qz^{-q})^{s-1}x^{q-1}y^{-q}(x^qy^{-q})^{s-1}X^{-t+1}
\end{split}
\end{equation}
Since $X^{-1}>1$, $y^qz^{-q}>1$ $z^{-1}>1$ and $y^{-1}>1$ then $y^qx^{-q}X^{-t}>1$. Hence $Y^ty^qx^{-q}X^{-t}>1$ and 
$Y^ty^qx^{-q}X^{-t+1}>1$ if $t>1$. If $t=1$, then relation $r''_1$ become
$$r'''_1=(Yy^qx^{-q}X^{-1})^{l-1}Yy^qx^{-q}Zz^qx^{-q}X^{-1}(Zz^qx^{-q}X^{-1})^{l-1}=1$$ and we have 

$y^qx^{-q}Z=y^qx^{-q}x^qz^{-q}(x^qz^{-q})^{s-1}z(y^qz^{-q})^{s}=y^qz^{-q}(x^qz^{-q})^{s-2}x^qz^{-q+1}(y^qz^{-q})^{s}>1.$
We have also by the first part of subcase (1) we have that 
 $X^tx^qz^{-q}Z^{-t}<1$ which implies $Z^tz^qx^{-q}X^{-t}=(X^tx^qz^{-q}Z^{-t})^{-1}>1$. This will gives a contradiction by the 
 relations $r''_1$ and $r'''_1$.
\end{itemize}

\end{enumerate}

 \textbf{Second case:} $q>0$, $s>0$, $t<0$, $l>0$.

We have $X=(y^qx^{-q})^{s}x(z^qx^{-q})^{s}=(y^qx^{-q})^{s-1}y^qx^{-q+1}(z^qx^{-q})^{s}$. Since 
$y<1$, $z<1$ and $x^{-1}<1$ then $X<1$. 

Similarly as the first case we have three subcases:

\begin{enumerate}
 \item If $Y<1$, then since $ZYX=1$ and $X<1$ then $Z>1$. Since 
 $Y=(z^qy^{-q})^{s}y(x^qy^{-q})^{s}=(z^qy^{-q})^{s}yxx^{q-1}y^{-q}(x^qy^{-q})^{s-1}
 =(z^qy^{-q})^{s}z^{-1}x^{q-1}y^{-q}(x^qy^{-q})^{s-1}$
then $z^qy^{-q}<1$ which implies $y^qz^{-q}>1$. By the relation 
$$r'_3=(X^tx^qz^{-q}Z^{-t})^{l}Z(Y^ty^qz^{-q}Z^{-t})^{l}=1$$ we have a contradiction.
\item If $Z<1$ then since $ZYX=1$ and $X<1$ then $Y>1$. 
 Since $Z=(x^qz^{-q})^{s}z(y^qz^{-q})^{s}=(x^qz^{-q})^{s-1}x^qz^{-q+1}(y^qz^{-q})^{s}$ then
 $y^qz^{-q}<1$ which implies $z^qy^{-q}>1$. By relation 
 $$r'_2=(Z^tz^qy^{-q}Y^{-t})^{l}Y(X^tx^qy^{-q}Y^{-t})^{l}=1$$ we have a contradiction.
\item The only subcase which is left to check is if $Y>1$ and $Z>1$. In this subcase the relation
$$r'_1=(Y^ty^qx^{-q}X^{-t})^{l}X(Z^tz^qx^{-q}X^{-t})^{l}=1$$ gives a contradiction.
 
\end{enumerate}

 \textbf{Third case:} $q<0$, $s>0$, $t<0$ and $l>0$.

We have $X=(y^qx^{-q})^{s}x(z^qx^{-q})^{s}>1$. Similarly as the first case we have three subcases:
\begin{enumerate}
 \item If $Y>1$, then since $ZYX=1$ and $X>1$ then $Z<1$. Since $Y=(z^qy^{-q})^{s}y(x^qy^{-q})^{s}$ then $z^qy^{-q}>1$ which
 implies that $y^qz^{-q}<1$. Therefore by the relation $$r'_3=(X^tx^qz^{-q}Z^{-t})^{l}Z(Y^ty^qz^{-q}Z^{-t})^{l}=1$$
 we have a contradiction.
 \item If $Z>1$, then since $ZYX=1$ and $X>1$ then $Y<1$. Since $Z=(x^qz^{-q})^{s}z(y^qz^{-q})^{s}$ then $y^qz^{-q}>1$ which
 implies that $z^qy^{-q}<1$. Therefore we have a contradiction by the relation
 $$r'_2=(Z^tz^qy^{-q}Y^{-t})^{l}Y(X^tx^qy^{-q}Y^{-t})^{l}=1$$ 
 \item The only subcase which is left to check is if $Y<1$ and $Z<1$. In this subcase the relation
$$r'_1=(Y^ty^qx^{-q}X^{-t})^{l}X(Z^tz^qx^{-q}X^{-t})^{l}=1$$ gives a contradiction.
\end{enumerate}

 \textbf{Fourth case:} $q<0$, $s<0$, $t<0$ and $l>0$.

We have $X=(x^qy^{-q})^{-s}x(x^qz^{-q})^{-s}=(x^qy^{-q})^{-s}x^{q+1}z^{-q}(x^qz^{-q})^{-s-1}<1$. 
Similarly as the first case we have three subcases:

\begin{enumerate}
 \item If $Y<1$, then since $ZYX=1$ and $X<1$ then $Z>1$. Since 
 $Y=(y^qz^{-q})^{-s}y(y^qx^{-q})^{-s}=(y^qz^{-q})^{-s}y^{q+1}x^{-q}(y^qx^{-q})^{-s-1}$ then $y^qz^{-q}<1$, so $z^qy^{-q}>1$.
 We have
 \begin{equation} \label{eq2}
\begin{split}
 Z^tz^qy^{-q} &=Z^{t+1}(y^qz^{-q})^{-s}z^{-1}(x^qz^{-q})^{-s}z^qy^{-q}\\
 &=Z^{t+1}(y^qz^{-q})^{-s}z^{-1}(x^qz^{-q})^{-s-1}x^qz^{-q}z^qy^{-q}\\
 &=
 Z^{t+1}(y^qz^{-q})^{-s}z^{-1}x^qz^{-q}(x^qz^{-q})^{-s-2}x^qy^{-q}\\
 &=
 Z^{t+1}(y^qz^{-q})^{-s}z^{-1}x^{-1}x^{q+1}z^{-q}(x^qz^{-q})^{-s-2}x^qy^{-q}\\
 &=
 Z^{t+1}(y^qz^{-q})^{-s}yx^{q+1}z^{-q}(x^qz^{-q})^{-s-2}x^qy^{-q}>1.
\end{split} 
 \end{equation}
 We have also that 
 \begin{equation} \label{eq3} 
\begin{split}
Z^tz^qx^{-q} &=Z^{t+1}(y^qz^{-q})^{-s}z^{-1}(x^qz^{-q})^{-s}z^qx^{-q}\\
 &=Z^{t+1}(y^qz^{-q})^{-s}z^{-1}(x^qz^{-q})^{-s-1}x^qz^{-q}z^qx^{-q}\\
 &=Z^{t+1}(y^qz^{-q})^{-s}z^{-1}(x^qz^{-q})^{-s-1}\\
 &=
 Z^{t+1}(y^qz^{-q})^{-s}z^{-1}x^qz^{-q}(x^qz^{-q})^{-s-2}\\
 &=
 Z^{t+1}(y^qz^{-q})^{-s}z^{-1}x^{-1}x^{q+1}z^{-q}(x^qz^{-q})^{-s-2}\\
 &=
 Z^{t+1}(y^qz^{-q})^{-s}yx^{q+1}z^{-q}(x^qz^{-q})^{-s-2}<1,\; \; {\rm if} \; s<-1.
\end{split}
\end{equation}
Therefore $x^qz^{-q}Z^{-t}>1$, and $X^tx^qz^{-q}Z^{-t}>1$ if $s<-1$. If $s=-1$ then
\begin{equation} \label{eq3} 
\begin{split}
Z^tz^qx^{-q}X^{-t} 
 &=Z^{t+1}(y^qz^{-q})z^{-1}XX^{-t-1}\\
 &=Z^{t+1}(y^qz^{-q})z^{-1}x^qy^{-q}xx^qz^{-q}X^{-t-1}\\
 &=Z^{t+1}(y^qz^{-q})z^{-1}x^{-1}x^{q+1}y^{-q}x^{q+1}z^{-q}X^{-t-1}\\
 &=Z^{t+1}(y^qz^{-q})yx^{q+1}y^{-q}x^{q+1}z^{-q}X^{-t-1}<1
\end{split}
\end{equation}
Hence $X^tx^qz^{-q}Z^{-t}>1$ if $s=-1$. 
Therefore by $r'_3$ we have a contradiction.
\item If $Z<1$ then since $ZYX=1$ and $X<1$ then $Y>1$. Since
\begin{equation}
\begin{split}
Z &=(z^qx^{-q})^{-s}z(z^qy^{-q})^{-s}\\
 &=(z^qx^{-q})^{-s-1}z^qx^{-q}z(z^qy^{-q})^{-s}\\
&=(z^qx^{-q})^{-s-1}z^qx^{-q-1}xz(z^qy^{-q})^{-s}\\
&=(z^qx^{-q})^{-s-1}z^qx^{-q-1}y^{-1}(z^qy^{-q})^{-s} 
\end{split}
\end{equation}
then $z^qy^{-q}<1$.

We have,
\begin{equation} \label{eq4} 
\begin{split}
Y^ty^qz^{-q}
&= Y^{t+1}(x^qy^{-q})^{-s}y^{-1}(z^qy^{-q})^{-s}y^qz^{-q}\\
&= Y^{t+1}(x^qy^{-q})^{-s-1}x^qy^{-q-1}(z^qy^{-q})^{-s-1}<1.
\end{split}
\end{equation}

Therefore $Z^tz^qy^{-q}Y^{-t}>1$. We have also that
\begin{equation} \label{eq5} 
\begin{split}
X^tx^qy^{-q}
&= X^{t+1}(z^qx^{-q})^{-s}x^{-1}(y^qx^{-q})^{-s-1}\\
&=X^{t+1}(z^qx^{-q})^{-s}z^qx^{-q-1}(y^qx^{-q})^{-s-1}>1
\end{split}
\end{equation}
Thus $X^tx^qy^{-q}Y^{-t}>1$. Therefore by $r'_2$ we have a contradiction.

\item The only subcase which is left to show is if $Y>1$ and $Z>1$. We have two subsubcases
\begin{itemize}
 \item If $z^qy^{-q}>1$ then by subcase (2) $Y^ty^qx^{-q}X^{-t}=(X^tx^qy^{-q}Y^{-t})^{-1}<1$ and by subcase (1) 
 $Z^tz^qx^{-q}X^{-t}=(X^tx^qz^{-q}Z^{-t})^{-1}<1$. Therefore we have a contradiction by $r'_1$.
 \item If $z^qy^{-q}<1$ then by subcase (2) we have $Y^ty^qx^{-q}X^{-t}=(X^tx^qy^{-q}Y^{-t})^{-1}<1$ and we have
 $Z^tz^qx^{-q}X^{-t}=Z^tz^qy^{-q}yqx^{-q}X^{-t}$ and by subcase (2) $yqx^{-q}X^{-t}<1$, so 
 $Z^tz^qx^{-q}X^{-t}=Z^tz^qy^{-q}yqx^{-q}X^{-t}<1$. Therefore we have a contradiction by $r'_1$.
\end{itemize}
\end{enumerate}

\textbf{ Fith case:} If $q<0$, $s>0$, $t>0$ and $l>0$.

We have $X=(y^qx^{-q})^{s}x(z^qx^{-q})^{s}>1$. 
Similarly as the first case we have three subcases:
\begin{enumerate}
 \item If $Y>1$ then since $ZYX=1$ and $X>1$ then $Z<1$. Since $Y=(z^qy^{-q})^{s}y(x^qy^{-q})^{s}>1$, so $z^qy^{-q}>1$.
 we have,
 \begin{equation}
 \begin{split}
 X^tx^qz^{-q} &= X^{t-1}(y^qx^{-q})^{s}x(z^qx^{-q})^{s}x^qz^{-q}\\
 &= X^{t-1}(y^qx^{-q})^{s}x(z^qx^{-q})^{s-1}z^qx^{-q}x^qz^{-q}\\
 &= X^{t-1}(y^qx^{-q})^{s}x(z^qx^{-q})^{s-1}>1
\end{split}
\end{equation}
Therefore $X^tx^qz^{-q}Z^{-t}>1$ and $X^tx^qz^{-q}Z^{-t+1}>1$. We have also that,
\begin{equation}
 \begin{split}
 Z^tz^qy^{-q}=Z^{t-1}(x^qz^{-q})^{s}z(y^qz^{-q})^{s-1}<1
\end{split}
\end{equation}
Hence $y^qz^{-q}Z^{-t}>1$ and $Y^ty^qz^{-q}Z^{-t}>1$, so we have a contradiction by $r''_3$.
\item If $Z>1$ then since $ZYX=1$ and $X>1$ then $Y<1$. Since $Z=(x^qz^{-q})^{s}z(y^qz^{-q})^{s}$, $y^qz^{-q}>1$. We have,
\begin{equation}
 \begin{split}
 Y^ty^qz^{-q} &= Y^{t-1}(z^qy^{-q})^{s}y(x^qy^{-q})^{s}y^qz^{-q}\\
 &= Y^{t-1}(z^qy^{-q})^{s}y(x^qy^{-q})^{s-1}x^qy^{-q}y^qz^{-q}\\
 &= Y^{t-1}(z^qy^{-q})^{s}y(x^qy^{-q})^{s-1}x^qz^{-q}<1
\end{split}
\end{equation}
Therefore $Y^ty^qz^{-q}Z^{-t}<1$, so $Z^tz^qy^{-q}Y^{-t}>1$ and $Z^tz^qy^{-q}Y^{-t+1}>1$ for $t>1$.
For $t=1$, the relation $r''_2$ become 
$$r'''_2=(Zz^qy^{-q}Y^{-1})^{l-1}Zz^qy^{-q}Xx^qy^{-q}Y^{-1}(Xx^qy^{-q}Y^{-1})^{l-1}=1$$ and we have
$$z^qy^{-q}X=z^qy^{-q}y^qx^{-q}(y^qx^{-q})^{s-1}x(z^qx^{-q})^{s}=z^qx^{-q}(y^qx^{-q})^{s-2}y^qx^{-q+1}(z^qx^{-q})^{s}>1.$$
We have also that,
\begin{equation}
 \begin{split}
 Y^ty^qx^{-q} &= Y^{t-1}(z^qy^{-q})^{s}y(x^qy^{-q})^{s}y^qx^{-q}\\
 &= Y^{t-1}(z^qy^{-q})^{s}y(x^qy^{-q})^{s-1}x^qy^{-q}y^qx^{-q}\\
 &= Y^{t-1}(z^qy^{-q})^{s}y(x^qy^{-q})^{s-1}<1
\end{split}
\end{equation}
Hence $Y^ty^qx^{-q}X^{-t}<1$ and $X^tx^qy^{-q}Y^{-t}>1$, so we have a contradiction by $r''_2$ and $r'''_2$.
\item The only subcase which is left to show is when $Y<1$ and $Z<1$. We have two subsubcases
\begin{itemize}
 \item If $z^qy^{-q}>1$ then by subcase (1) $Z^tz^qx^{-q}X^{-t}<1$. We have,
 \begin{equation}
 \begin{split}
 X^tx^qy^{-q} &= X^{t-1}(y^qx^{-q})^{s}x(z^qx^{-q})^{s}x^qy^{-q}\\
 &= X^{t-1}(y^qx^{-q})^{s}x(z^qx^{-q})^{s-1}z^qx^{-q}x^qy^{-q}\\
 &= X^{t-1}(y^qx^{-q})^{s}x(z^qx^{-q})^{s-1}z^qy^{-q}>1
\end{split}
\end{equation}
Therefore $X^tx^qy^{-q}Y^{-t}>1$, so $Y^ty^qx^{-q}X^{-t}<1$ and $Y^ty^qx^{-q}X^{-t+1}<1$ for $t>1$. Hence we have a 
contradiction by $r''_1$ if $t>1$. If $t=1$ use a similar argument as the previous case to conclude.
\item If $z^qy^{-q}<1$ then $y^qz^{-q}>1$. We have $Y^ty^qx^{-q}X^{-t}<1$ and $Y^ty^qx^{-q}X^{-t+1}<1$ by subcase (2). We also
have $Z^tz^qx^{-q}X^{-t}<1$ by subcase (1). Therefore we have a contradiction by $r'_1$.
\end{itemize}
\end{enumerate}

 \textbf{Sixth case:} $q<0$, $s>0$, $t<0$ and $l<0$.

We have $X=(y^qx^{-q})^{s}x(z^qx^{-q})^{s}>1$. 
Similarly as the first case we have three subcases:
\begin{enumerate}
 \item If $Y>1$ then since $ZYX=1$ then $Z<1$. Since $Y=(z^qy^{-q})^{s}y(x^qy^{-q})^{s}>1$, $z^qy^{-q}>1$. The relation
 $r'_3$ become
 $r'''_3=
 (X^tx^qz^{-q}Z^{-t})^{l}Z^{t+1}z^qy^{-q}Y^{-t}(Y^ty^qz^{-q}Z^{-t})^{l+1}=1$ and this implies a contradiction.
 \item If $Z>1$ then since $ZYX=1$ then $Y<1$. Since $Z=(x^qz^{-q})^{s}z(y^qz^{-q})^{s}$, $y^qz^{-q}>1$. The relation
 $r'_2$ become
 $r'''_2=
 (Z^tz^qy^{-q}Y^{-t})^{l}Y^{t+1}y^qx^{-q}X^{-t}(X^tx^qy^{-q}Y^{-t})^{l+1} =1$, which implies a contradiction.
 \item The only subcase which is left to show is when $Y<1$ and $Z<1$ and in this 
 subcase we have a contradiction by the following:
 $r'''_1=
 (Y^ty^qx^{-q}X^{-t})^{l}X^{t+1}x^qz^{-q}Z^{-t}(Z^tz^qx^{-q}X^{-t})^{l+1}=1$.
\end{enumerate}

\textbf{ Seventh case:} $q>0$, $s<0$, $t<0$ and $l>0$.

We have $X=(y^qx^{-q})^{s}x(z^qx^{-q})^{s}>1$. 
Similarly as the first case we have three subcases:
\begin{enumerate}
 \item If $Y>1$ then since $ZYX=1$ then $Z<1$. Since $Y=(z^qy^{-q})^{s}y(x^qy^{-q})^{s}>1$, $z^qy^{-q}<1$. 
 We have 
\begin{equation} \label{eq1}
\begin{split}
x^qz^{-q}Z^{-t} &= x^qz^{-q}(x^qz^{-q})^{s}z(y^qz^{-q})^{s}Z^{-t-1}\\
 &= x^qz^{-q}(z^qx^{-q})^{-s}z(y^qz^{-q})^{s}Z^{-t-1}\\
 &=x^qz^{-q}z^qx^{-q}(z^qx^{-q})^{-s-1}z(y^qz^{-q})^{s}Z^{-t-1}\\
 &=(z^qx^{-q})^{-s-1}z(y^qz^{-q})^{s}Z^{-t-1}<1
\end{split}
\end{equation}
Therefore $X^tx^qz^{-q}Z^{-t}<1$.
We have also that,
\begin{equation} \label{eq1}
\begin{split}
y^qz^{-q}Z^{-t} &= y^qz^{-q}(x^qz^{-q})^{s}z(y^qz^{-q})^{s}Z^{-t-1}\\
 &= y^qz^{-q}(z^qx^{-q})^{-s}z(y^qz^{-q})^{s}Z^{-t-1}\\
 &=y^qz^{-q}z^qx^{-q}(z^qx^{-q})^{-s-1}z(y^qz^{-q})^{s}Z^{-t-1}\\
 &=y^qx^{-q}(z^qx^{-q})^{-s-1}z(y^qz^{-q})^{s}Z^{-t-1}<1
\end{split}
\end{equation}
Therefore $Y^ty^qz^{-q}Z^{-t}<1$. Thus we have a contradiction by $r'_3$.
\item If $Z>1$ then since $ZYX=1$ then $Y<1$. Since $Z=(x^qz^{-q})^{s}z(y^qz^{-q})^{s}$, $y^qz^{-q}<1$. We have,
 \begin{equation}
 \begin{split}
 z^qy^{-q}Y^{-t} &= z^qy^{-q}(z^qy^{-q})^{s}y(x^qy^{-q})^{s}Y^{-t-1}\\
 &= (z^qy^{-q})^{s+1}y(x^qy^{-q})^{s-1}Y^{-t-1}<1
\end{split}
\end{equation}
Hence $Z^tz^qy^{-q}Y^{-t}<1$. We also have,
\begin{equation}
 \begin{split}
 X^tx^qy^{-q} &= X^{t+1}(z^qx^{-q})^{-s}x^{-1}(y^qx^{-q})^{-s}x^qy^{-q}\\
 &= X^{t+1}(z^qx^{-q})^{-s}x^{-1}(y^qx^{-q})^{-s-1}<1
\end{split}
\end{equation}
Therefore $X^tx^qy^{-q}Y^{-t}<1$, and we have a contradiction by $r'_2$.
\item The only subcase which is left to show is when $Y<1$ and $Z<1$ and in this subcase we have two subsubcases:
 \begin{itemize}
 \item If $z^qy^{-q}>1$ then by subcase (2) we have $X^tx^qy^{-q}<1$, so $y^qx^{-q}X^{-t}>1$. Therefore we have a contradiction
 by the following:
 $$r'_1=(Y^ty^qx^{-q}X^{-t})^{l}X(Z^tz^qy^{-q}y^qx^{-q}X^{-t})^{l}=1$$
 \item If $z^qy^{-q}<1$ then by subcase (1) $X^tx^qz^{-q}Z^{-t}<1$, so $Z^tz^qx^{-q}X^{-t}>1$. We have also by subcase (2) 
 $X^tx^qy^{-q}Y^{-t}<1$, so $Y^ty^qx^{-q}X^{-t}>1$, and we have a contradiction by $r'_1$.
\end{itemize}
\end{enumerate}

\textbf{ Eighth case:} $q>0$, $s>0$, $t<0$ and $l<0$.

We have $X=(y^qx^{-q})^{s}x(z^qx^{-q})^{s}=(y^qx^{-q})^{s-1}y^qx^{-q+1}(z^qx^{-q})^{s}<1$.

We discuss the signs of $Y$ and $Z$. We have three subcases
\begin{enumerate}
 \item If $Y<1$, then since $ZYX=1$ and $X<1$ then $Z>1$. Since 
 $Y=(z^qy^{-q})^{s}y(x^qy^{-q})^{s}=(z^qy^{-q})^{s}yxx^{q-1}y^{-q}(x^qy^{-q})^{s-1}
 =(z^qy^{-q})^{s}z^{-1}x^{q-1}y^{-q}(x^qy^{-q})^{s-1}$
then $z^qy^{-q}<1$, and $y^qz^{-q}>1$. 
We have a contradiction by the following:
$$r'_3=
 (X^tx^qz^{-q}Z^{-t})^{l}Z^{t+1}z^qy^{-q}Y^{-t}(Y^ty^qz^{-q}Z^{-t})^{l+1}=1$$
 \item If $Z<1$ then since $ZYX=1$ and $X<1$ then $Y>1$. 
 Since $Z=(x^qz^{-q})^{s}z(y^qz^{-q})^{s}=(x^qz^{-q})^{s-1}x^qz^{-q+1}(y^qz^{-q})^{s}$ then
 $y^qz^{-q}<1$, so $z^qy^{-q}>1$. We have a contradiction by the following:
 $$r'_2=
 (Z^tz^qy^{-q}Y^{-t})^{l}Y^{t+1}y^qx^{-q}X^{-t}(X^tx^qy^{-q}Y^{-t})^{l+1} =1$$
\item The only subcase which is left to check is if $Y>1$ and $Z>1$. In this subcase we have a contradiction by the following:
 $$r'_1=
 (Y^ty^qx^{-q}X^{-t})^{l}X^{t+1}x^qz^{-q}Z^{-t}(Z^tz^qx^{-q}X^{-t})^{l+1}=1$$.
\end{enumerate}
This complete the proof of Theorem \ref {thm: main result2}.

\end{proof}

 \begin{remark}
{\rm 
\begin{enumerate}
 \item Since the knot $5_1$ corresponds to the two-bridge knot $K_{[-2,2,-2,2]}$ and, by Gordon and 
Lidman [GL], the fundamental group of the $n$-fold cyclic branched cover of $5_1$ is left-orderable for $n\geq 4$, 
then another question we can also ask is the 
following: Is the fundamental group of the n-fold cyclic branched cover of $K_{[-2q,2s,-2t,2l]}$
left-orderable for $n\geq 4$?
\item In general, it is not true that for every 2-bridge knot $K$ the fundamental group of the 3-fold cyclic branched cover of
$K$ is not left-orderable,
as an example we have the 2-bridge knot $K_{[6,-3]}$ which can also be written as the genus three 2-bridge knot 
$K_{[-2,2,-2,2,-2,-2]}$. By [GL] this knot has left-orderable fundamental group. The proof come from the fact that 
$\pi_1(\Sigma_3(K_{[6,-3]}))$ has a nontrivial representation into $PSL(2,\mathbb{R})$. Since $\Sigma_3(K_{[6,-3]})$ is an 
integer homology sphere, this representation lifts to a nontrivial representation into $\widetilde{SL}(2,\mathbb{R})$.
Since $\widetilde{SL}(2,\mathbb{R})$ is left-orderable then $\pi_1(\Sigma_3(K_{[6,-3]}))$ is also left-orderable.
\end{enumerate}}
\end{remark}
\section{L-spaces and genus 2 two-bridge knots}
In this section we will complete the proof of Theorem \ref{thm: main result4} by proving the following theorem.

\begin{theorem}\label{thm: main result1}
  The 3-fold cyclic branched cover of $K_{[-2q,2s,-2t,2l]}$ is an L-space, where $q$, $s$, $t$ and $l$ 
  $\in\mathbb{Z}\setminus\{0\}$.
\end{theorem}

In order to show that the 3-fold cyclic branched cover of any genus 2 two-bridge knot $K_{[-2q,2s,-2t,2l]}$ 
is an L-space we have to consider all cases for the signs of $q$, $s$, $t$ and $l$. We have sixteen such cases, 
but since the mirror image of $K_{[-2q,2s,-2t,2l]}$ is $K_{[2q,-2s,2t,-2l]}$ we need only 
deal with eight of them:
\begin{enumerate}
 \item $q>0$, $s>0$, $t>0$ and $l>0$;
 \item $q<0$, $s>0$, $t<0$ and $l>0$;
 \item $q>0$, $s<0$, $t>0$ and $l>0$;
 \item $q<0$, $s<0$, $t<0$ and $l>0$;
 \item $q>0$, $s>0$, $t<0$ and $l>0$;
 \item $q>0$, $s<0$, $t<0$ and $l<0$;
 \item $q>0$, $s<0$, $t<0$ and $l>0$;
 \item $q<0$, $s<0$, $t>0$ and $l>0$.
 
\end{enumerate}

By $A(t;\ast,\ast,\ast)$ we mean the link $A$ pictured in Figure 3. Let $\gamma \in \{\infty, 0\}$.
By $A(t;\gamma,\ast,\ast)$ we mean the link $A$ with the resolution $\gamma$ (cf. Figure 2) at
the leftmost $t$ half twists (see Figures 5 and 6 for the case $t>1$). 
By $A(t;\ast,\gamma,\ast)$ we mean the link $A$ with the resolution $\gamma$ at
the $t$ half twists on the middle (see Figures 7, 8, 9 and 10 for the case $t>1$)
and by $A(t;\ast,\ast,\gamma)$ we mean the link $A$ with the resolution $\gamma$ at
the $t$ half twists on the right (see Figures 11, 12, 13, 14, 15 and 16 for the case $t>1$). Similar notation will be used also
for $L(l;\ast,\ast,\ast)$ (cf. Figure 17).

\begin{figure}[H]
\centering
\def\svgwidth{0.50\columnwidth}
 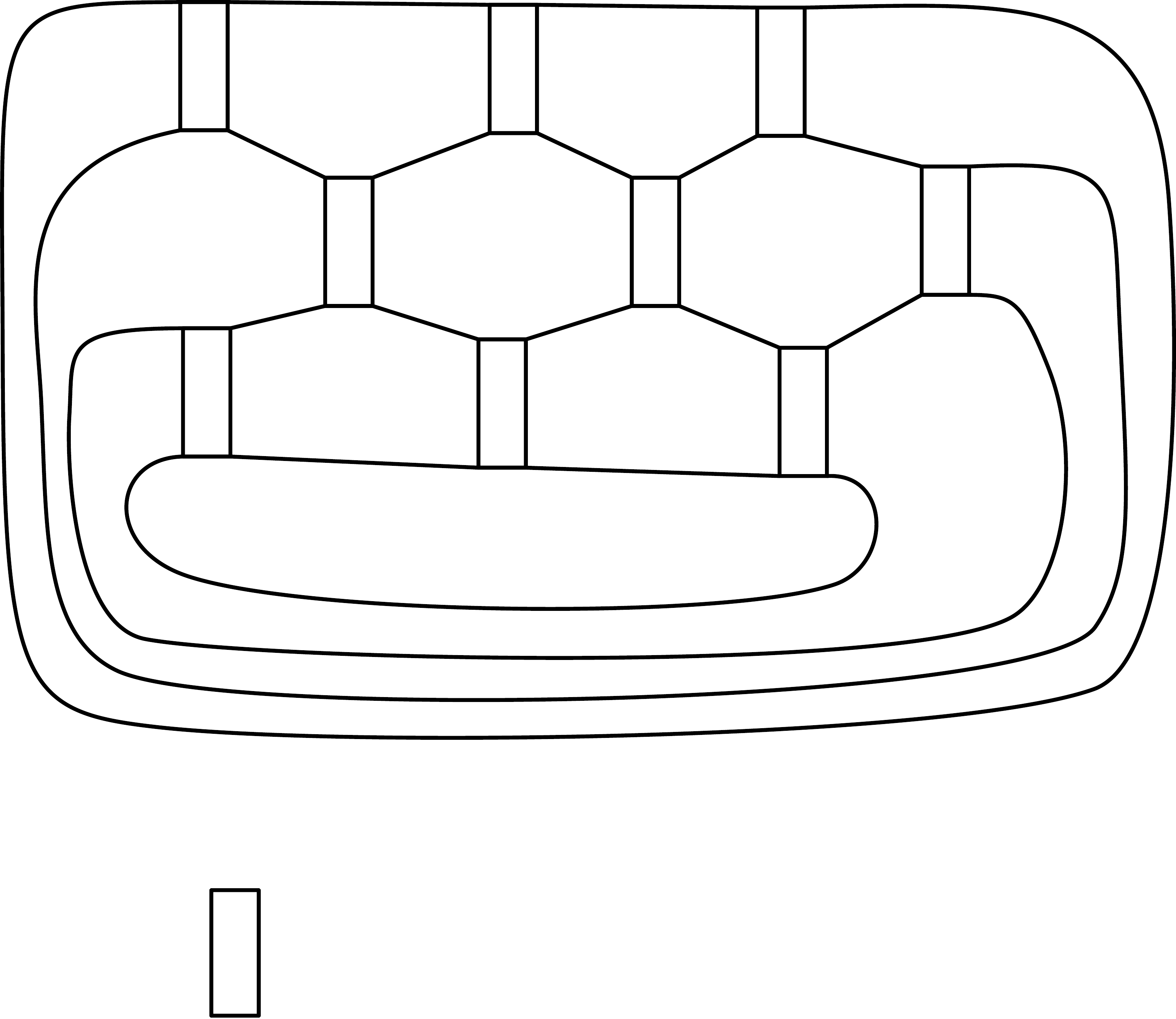
 \caption{The link $A$.}
\label{fig:The link A.}
\end{figure}
\begin{figure}[H]
\centering
\def\svgwidth{0.50\columnwidth}
 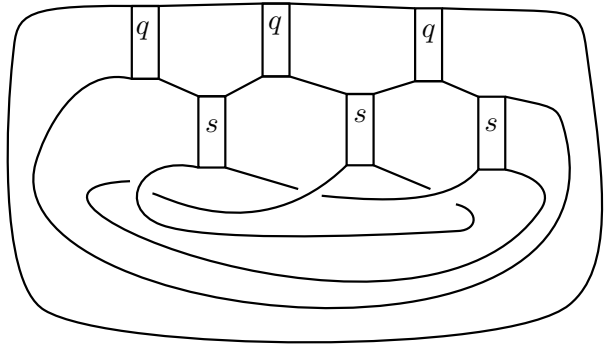
 \caption{The link $A$ for $t=1$}
\label{fig:The link A for t=1}
\end{figure}
The proof of Theorem \ref{thm: main result1} will be split into two parts. First, we will study the link
$A(t;\ast,\ast,\ast)$. That is, we study
when it is quasi-alternating or when $\Sigma_2( A(t;\ast,\ast,\ast))$ is an $L$-space. Second, we use the link
$A(t;\ast,\ast,\ast)$
to show that $\Sigma_2( L(l;\ast,\ast,\ast))$ is an $L$-space.
\subsection{The link $A(t;\ast,\ast,\ast)$}

We have 4 cases,

1) If $q>0$, $s>0$, $t>0$ 
\begin{figure}[H]
\centering
\def\svgwidth{0.50\columnwidth}
 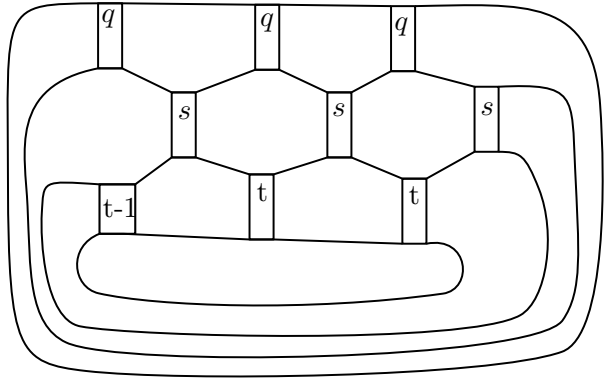
 \caption{The link $A(t;\infty,\ast,\ast)$.}
\label{fig:The link A.}
\end{figure}
\begin{figure}[H]
\centering
\def\svgwidth{0.50\columnwidth}
 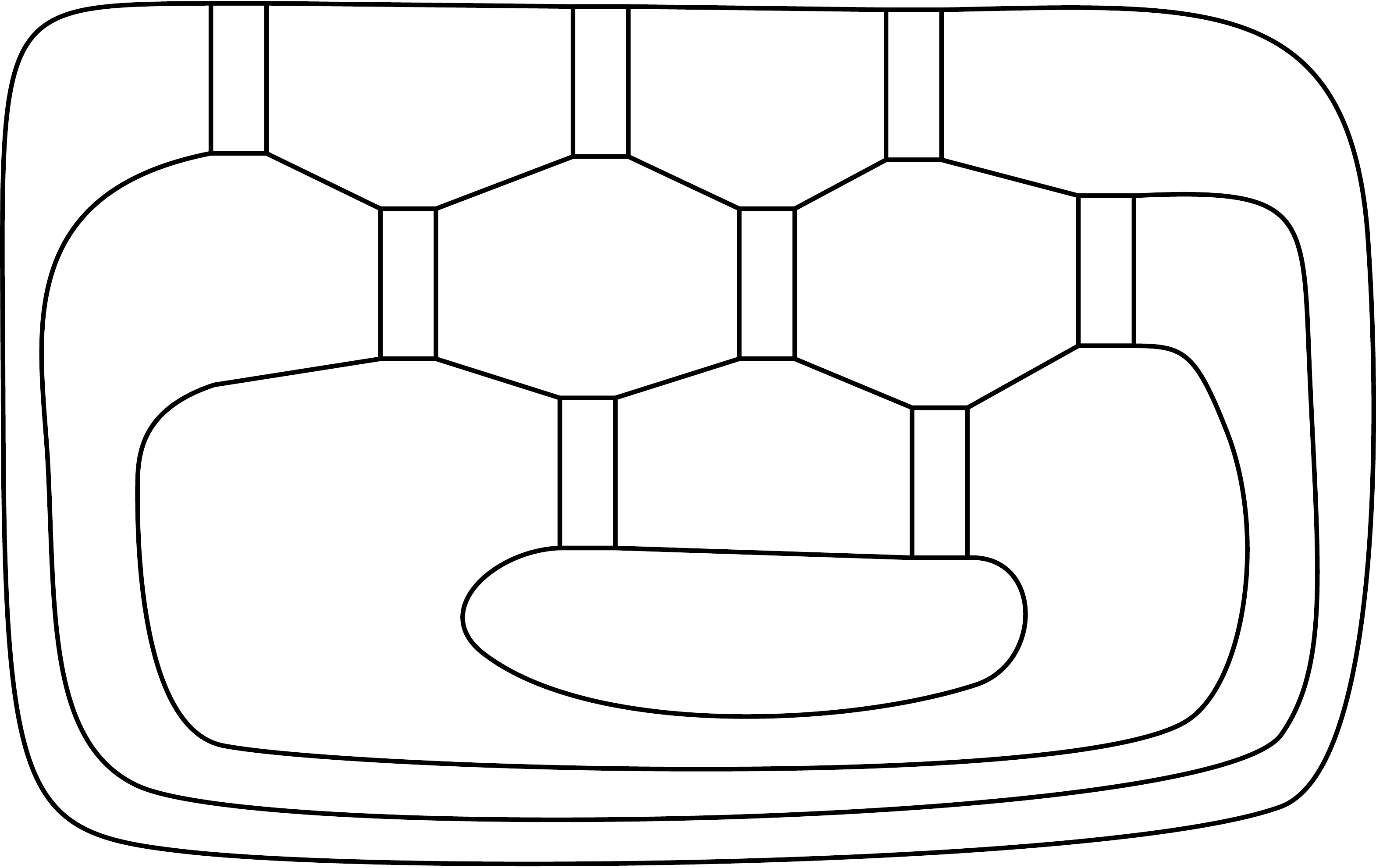
 \caption{The link $A(t;0,\ast,\ast)$.}
\label{fig:The link A.}
\end{figure}
\begin{figure}[H]
\centering
\def\svgwidth{0.50\columnwidth}
 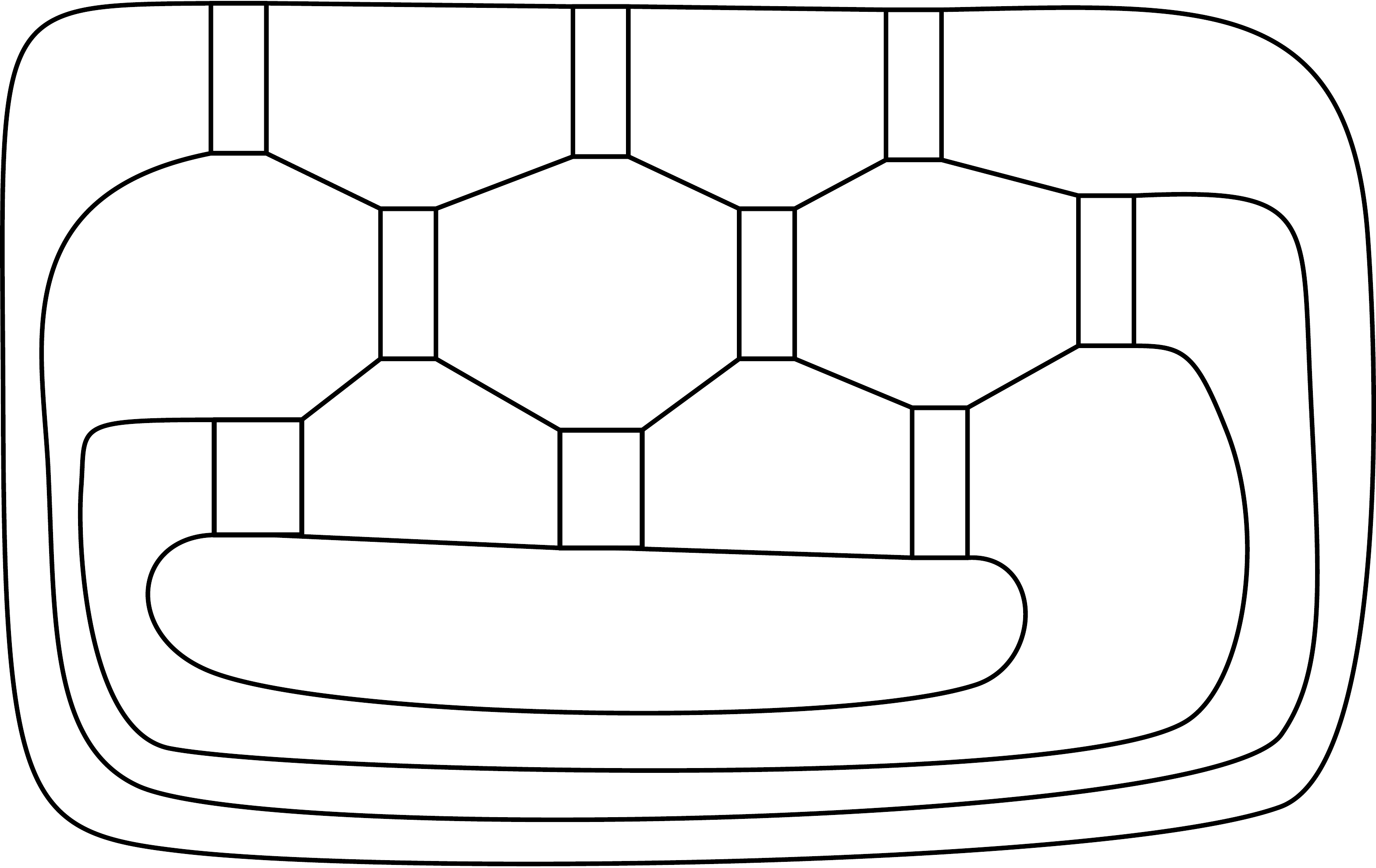
 \caption{The link $A(t;\infty,\infty,\ast)$.}
\label{fig:The link A.}
\end{figure}
\begin{figure}[H]
\centering
\def\svgwidth{0.50\columnwidth}
 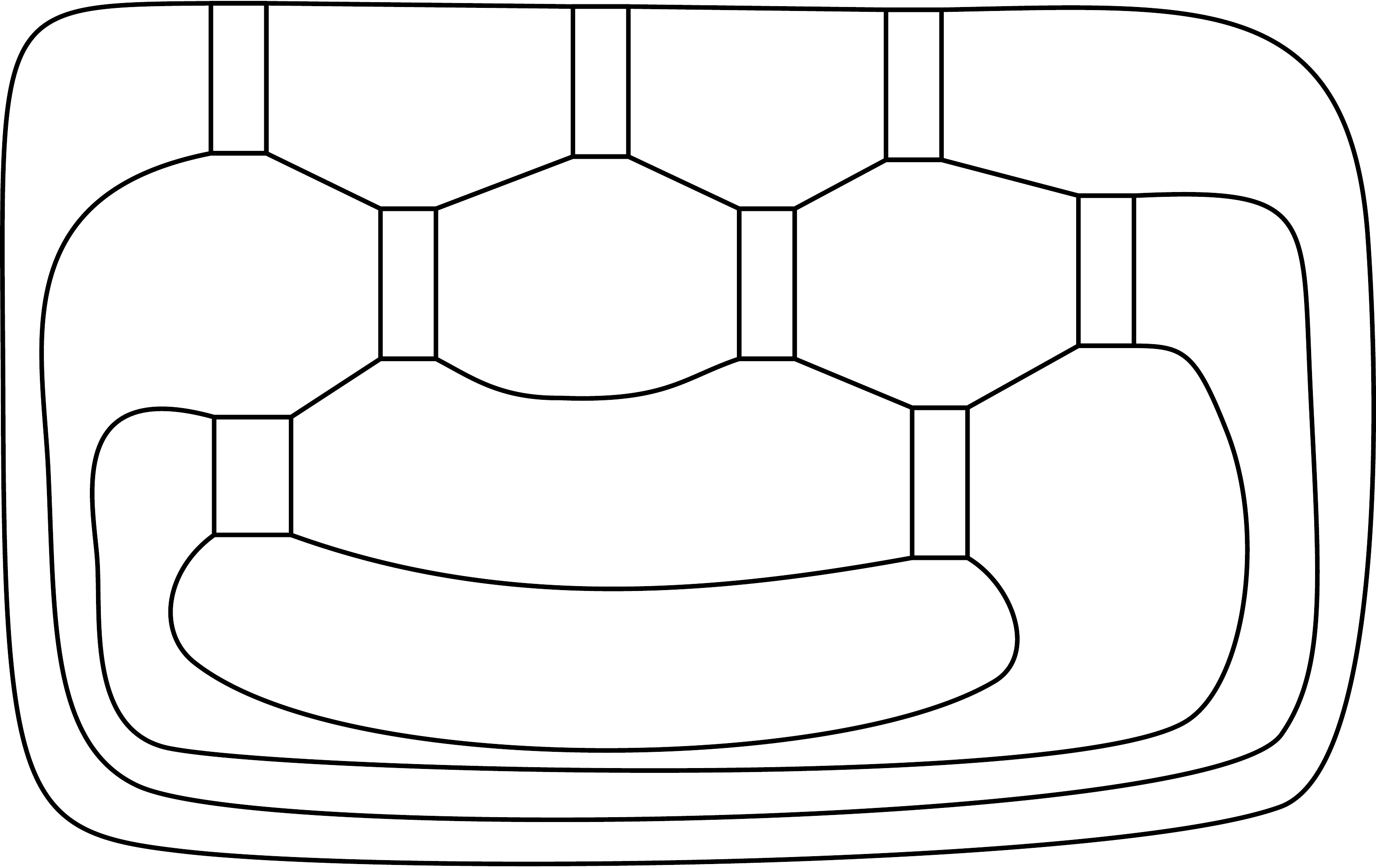
 \caption{The link $A(t;\infty,0,\ast)$.}
\label{fig:The link A.}
\end{figure}
\begin{figure}[H]
\centering
\def\svgwidth{0.50\columnwidth}
 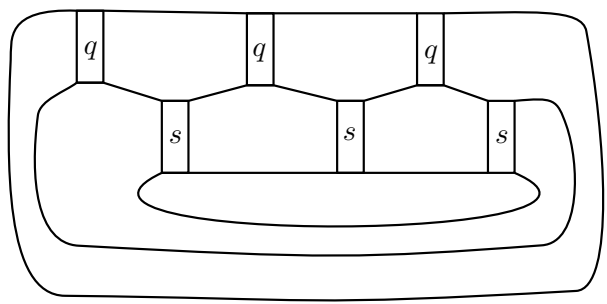
 \caption{The link $A(t;0,0,\ast)$.}
\label{fig:The link A.}
\end{figure}
\begin{figure}[H]
\centering
\def\svgwidth{0.50\columnwidth}
 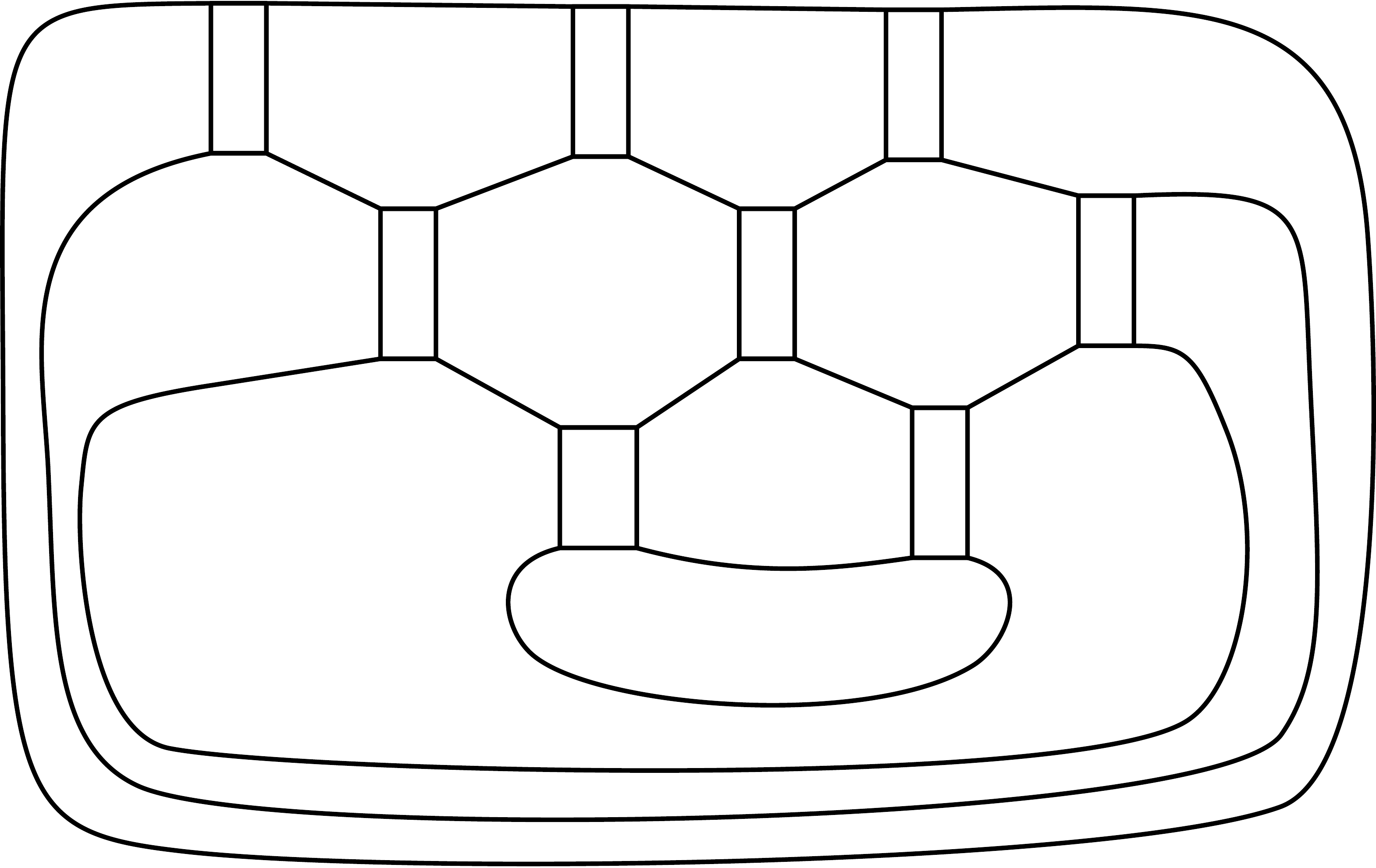
 \caption{The link $A(t;0,\infty,\ast)$.}
\label{fig:The link A.}
\end{figure}
\begin{figure}[H]
\centering
\def\svgwidth{0.50\columnwidth}
 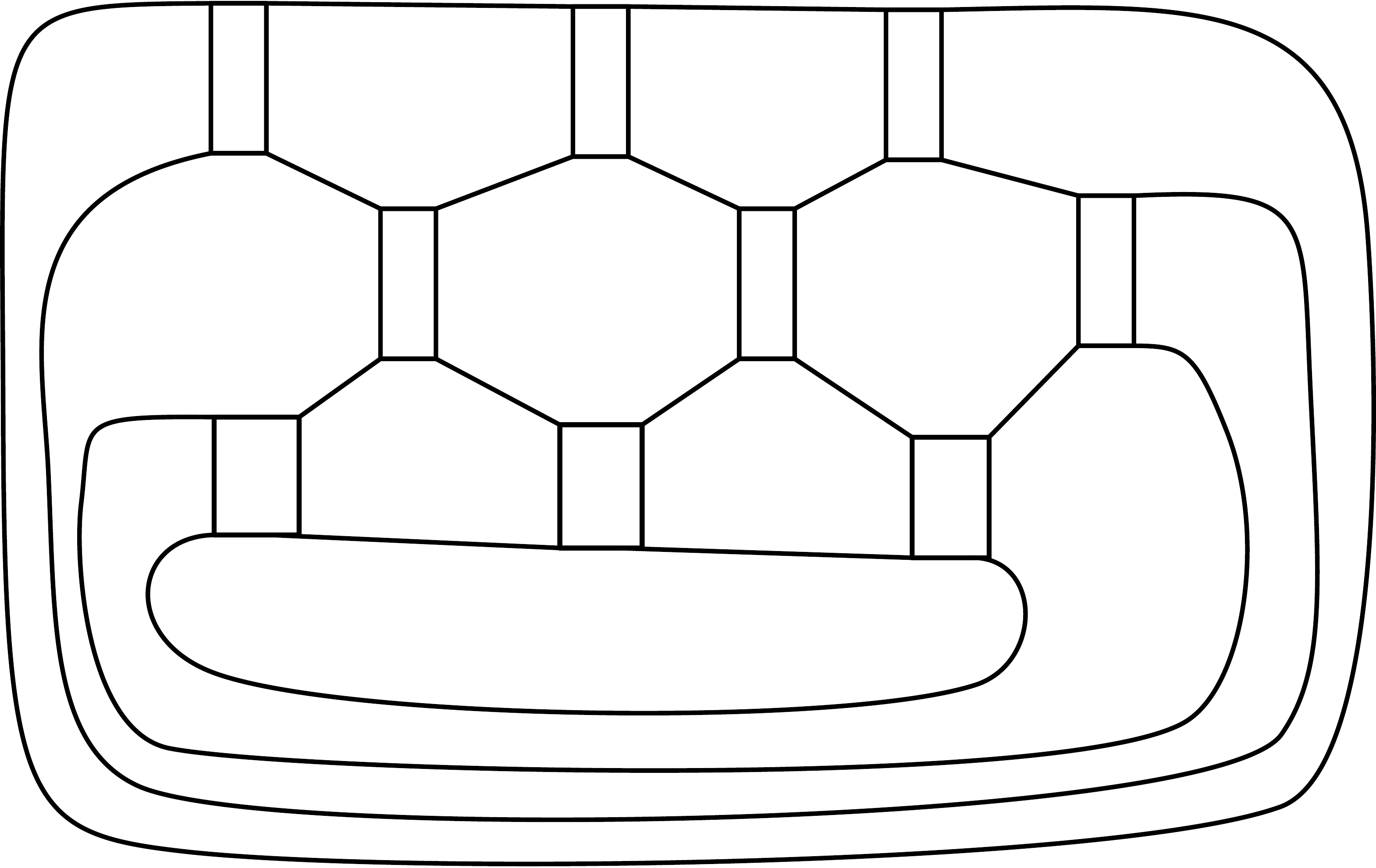
 \caption{The link $A(t;\infty,\infty,\infty)$.}
\label{fig:The link A.}
\end{figure}
\begin{figure}[H]
\centering
\def\svgwidth{0.50\columnwidth}
 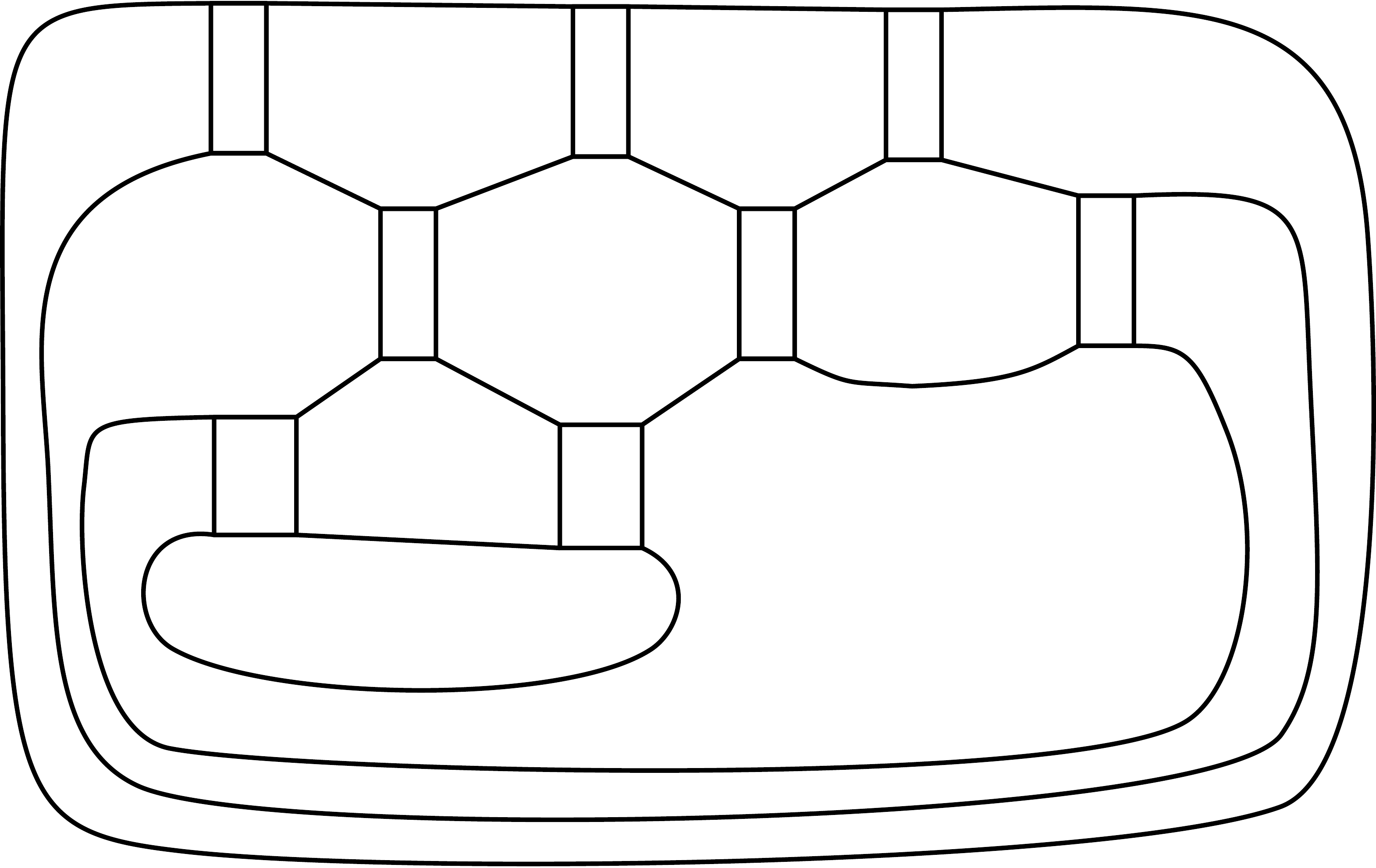
 \caption{The link $A(t;\infty,\infty,0)$.}
\label{fig:The link A.}
\end{figure}
\begin{figure}[H]
\centering
\def\svgwidth{0.50\columnwidth}
 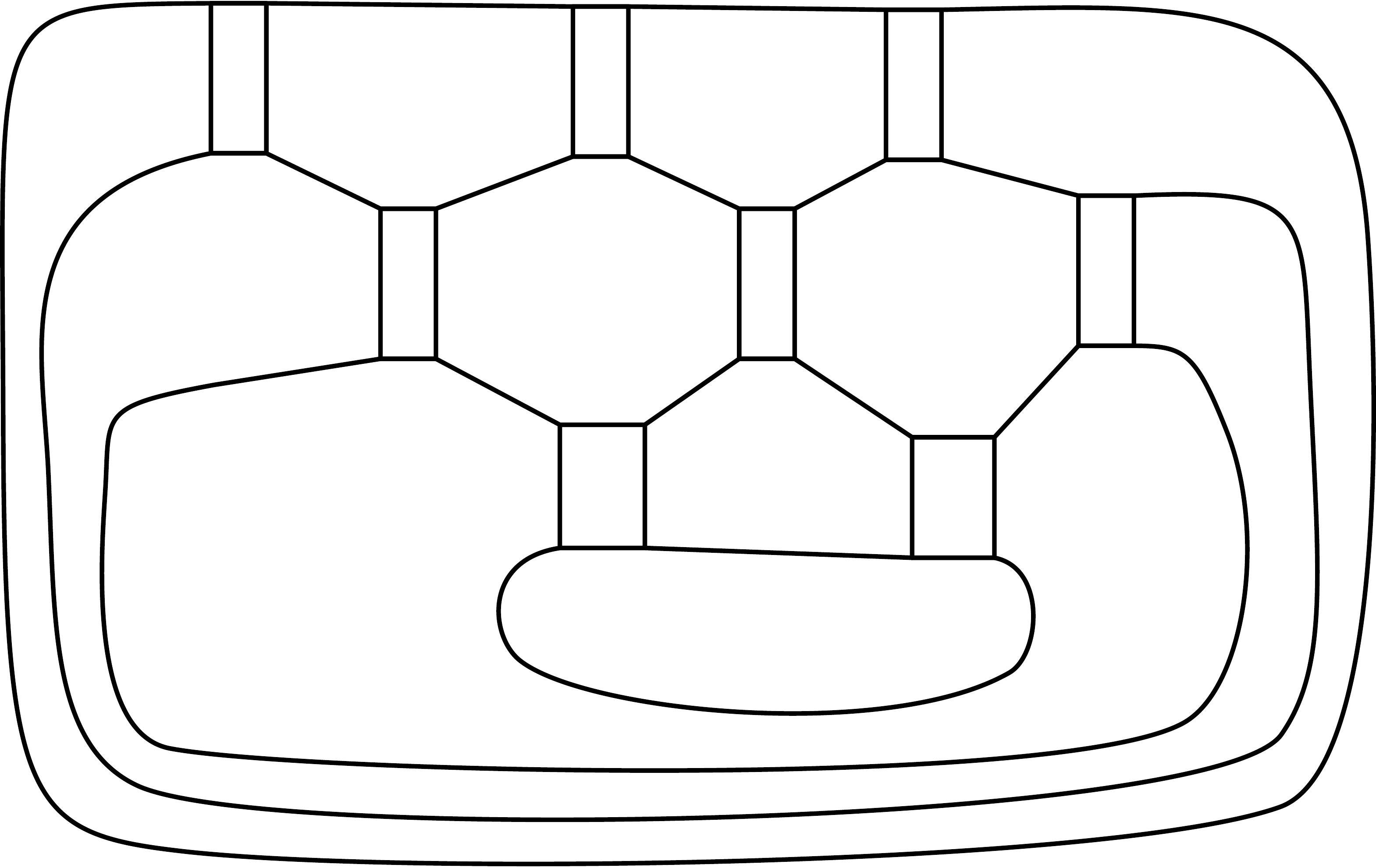
 \caption{The link $A(t;0,\infty,\infty)$.}
\label{fig:The link A.}
\end{figure}
\begin{figure}[H]
\centering
\def\svgwidth{0.50\columnwidth}
 \input{A12.eps_tex}
 \caption{The link $A(t;0,\infty,0)$.}
\label{fig:The link A.}
\end{figure}
\begin{figure}[H]
\centering
\def\svgwidth{0.50\columnwidth}
 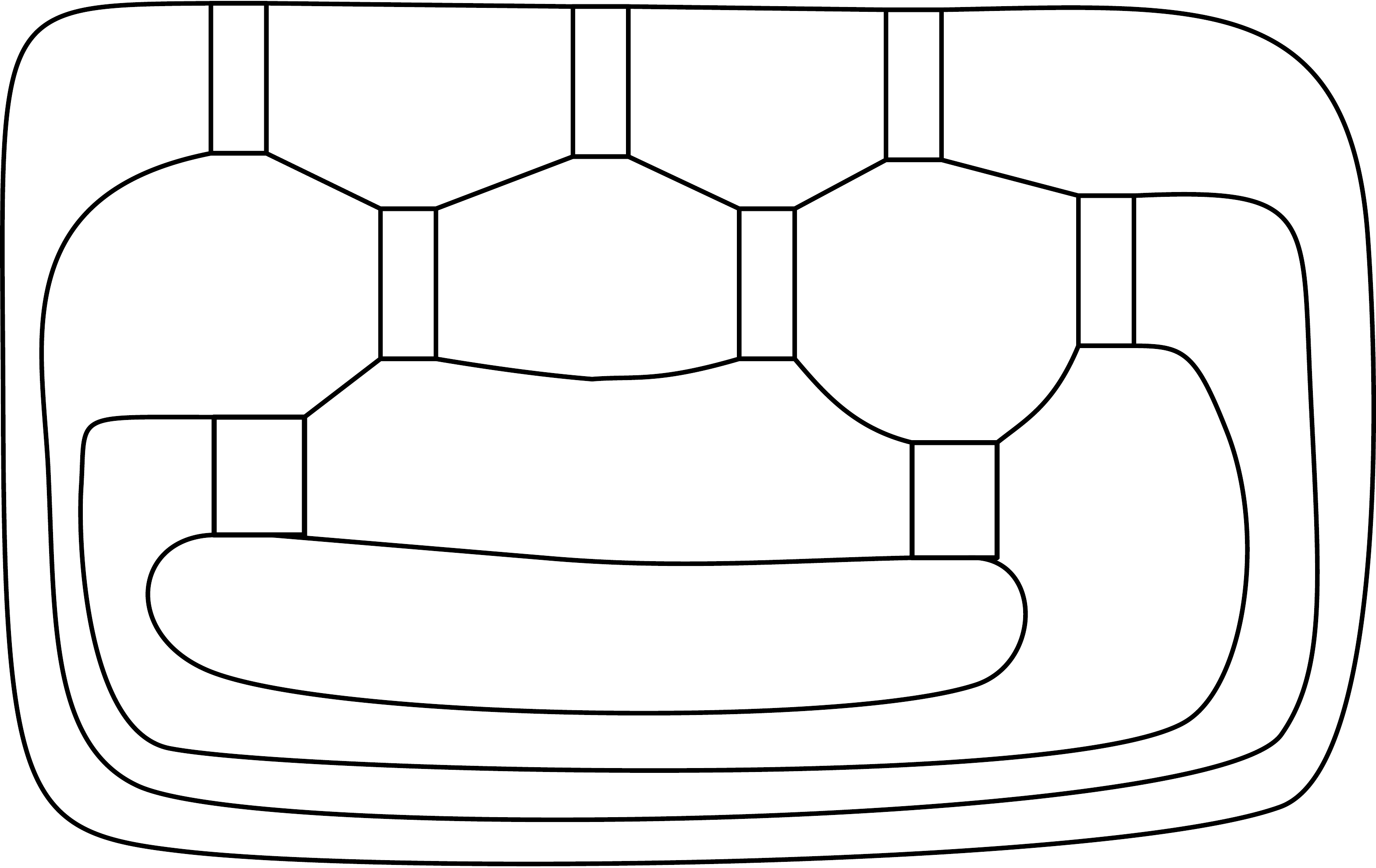
 \caption{The link $A(t;\infty,0,\infty)$.}
\label{fig:The link A.}
\end{figure}
\begin{figure}[H]
\centering
\def\svgwidth{0.50\columnwidth}
 \input{A12.eps_tex}
 \caption{The link $A(t;\infty,0,0)$.}
\label{fig:The link A.}
\end{figure}
The Goeritz matrix of $A(t;\ast,\ast,\ast)$ is 
\[
  \begin{bmatrix}
    -2I & I & O & \dots  & O & \dots \\
    I & -2I & I & O & \dots  & O & \dots \\
    \vdots & \vdots & \vdots & \ddots & \vdots & \vdots & \vdots \\
    O & O \dots   &    O & I & S & I& O \dots\\ 
     O & \dots  & O& O & I & -2I & I & O  \dots \\
     \vdots & \vdots & \vdots & & \vdots & \ddots & \vdots & \vdots & \vdots \\
     O & \dots  & O & O & \dots O &  I & -2I & I & F \\
    O & O & O & \dots  & \dots & O  & I & -2I & E \\
    Z & Z & \dots & \dots & Z & Z & Z & D & -3
\end{bmatrix}
\]

where \[
   S=
  \left[ {\begin{array}{ccc}
   2s-2 & -s & -s \\
   -s & 2s-2 & -s \\
   -s & -s & 2s-2
  \end{array} } \right],
  F=
  \left[ {\begin{array}{c}             
   0 \\
   0 \\
   0
  \end{array} } \right]  
\] 

\[
   E=
  \left[ {\begin{array}{c}
   1 \\
   1 \\
   1
  \end{array} } \right]
\]

\[
   O=
  \left[ {\begin{array}{ccc}
   0 & 0 & 0 \\
   0 & 0 & 0 \\
   0 & 0 & 0
  \end{array} } \right],
  I=
  \left[ {\begin{array}{ccc}
   1 & 0 & 0 \\
   0 & 1 & 0 \\
   0 & 0 & 1
  \end{array} } \right]
\] 
\[
   Z=
  \left[ {\begin{array}{ccc}
   0 & 0 & 0
  \end{array} } \right],
   D=
  \left[ {\begin{array}{ccc}
   1 & 1 & 1
  \end{array} } \right]
\]

a) By induction, assume $t=1$.

In the following table we give the determinant of the corresponding link, obtained through brute calculation.
\begin{table}[H]\label{Table: 2}
\centering
\begin{tabular}{ | m{3cm} | m{8cm} | } 
\hline
Link & Determinant ($s>1$)  \\ 
\hline
$A(t=1;\ast,\ast,\ast)$ & $3(-1-q+3qs)^2$ \\ 
\hline
 $A(t=1;0,\ast,\ast)$ & $2(-1-q+3qs)(-1+3qs)$  \\ 
\hline
$A(t=1;\infty,\ast,\ast)$ & $(-1-q+3qs)(-1-3q+3qs)$  \\ 
\hline
$A(t=1;0,\infty,\ast)$ & $(-1+3qs)(-1-2q+3qs)$   \\ 
\hline
 $A(t=1;0,0,\ast)$& $(-1+3qs)^2$   \\ 
\hline
\end{tabular}
\caption{}
\end{table}

\begin{lemma}\label{lemma: 5.2}
 The links $A(t=1;\infty,\ast,\ast)$, $A(t=1;0,\infty,\ast)$ and $A(t=1;0,0,\ast)$ are quasi-alternating for $s> 1$.
 Therefore, $A(t=1;0,\ast,\ast)$ and $A(t=1;\ast,\ast,\ast)$
 are also quasi-alternating for $s>1$.
\end{lemma}
\begin{proof}
 The fact that the links $A(t=1;\infty,\ast,\ast)$, $A(t=1;0,\infty,\ast)$ and 
 $A(t=1;0,0,\ast)$ are quasi-alternating for $s> 1$ was shown by Peters
 [P]. Since det $A(t=1;0,\ast,\ast)$ = det $A(t=1;0,\infty,\ast)$ + det $A(t=1;0,0,\ast)$ and 
 det $A(t=1;\ast,\ast,\ast)$ = det $A(t=1;0,\ast,\ast)$ + det $A(t=1;\infty,\ast,\ast)$ for $s> 1$ by Table 2,
 then $A(t=1;0,\ast,\ast)$ and $A(t=1;\ast,\ast,\ast)$ are
 quasi-alternating for $s> 1$.
\end{proof}

b) Assume $t>1$.

We have the following table which gives the determinant of the corresponding link.
\begin{table}[H]\label{Table: 3}
\centering
\begin{tabular}{ | m{3cm} | m{8cm} | } 
\hline
Link & Determinant ($t>1$) \\ 
\hline
$A(t;\ast,\ast,\ast)$ & $3(-t-q+3qst)^2$ \\ 
\hline
 $A(t;0,\ast,\ast)$ & $2(-t-q+3qst)(-1+3qs)$  \\ 
\hline
$A(t;\infty,\ast,\ast)$ & $(-t-q+3qst)(2-3q-6qs-3t+9qst)$  \\ 
\hline
$A(t;0,\infty,\ast)$ & $(-1+3qs)(1-2q-3qs-2t+6qst)$   \\ 
\hline
$A(t;\infty,0,\ast)$ & $(-1+3qs)(1-2q-3qs-2t+6qst)$   \\ 
\hline
 $A(t;0,0,\ast)$& $(-1+3qs)^2$   \\ 
\hline
$A(t;\infty,\infty,\ast)$ & $(1-q-3qs-t+3qst)(1-3q-3qs-3t+9qst)$  \\ 
\hline
$A(t;\infty,\infty,\infty)$ & $3(1-q-3qs-t+3qst)^2$  \\ 
\hline
$A(t;0,\infty,\infty)$ & $2(-1+3qs)(1-q-3qs-t+3qst)$  \\ 
\hline
\end{tabular}
\caption{}
\end{table}

\begin{lemma}\label{lemma: 5.3}
\begin{enumerate}
 \item $A(t;0,0,\ast)$ = $A(t;\infty,0,0)$ = $A(t;0,\infty,0)$, 
 \item $A(t;\infty,0,\infty)$ = $A(t;\infty,\infty,0)$ = $A(t;0,\infty,\infty)$, 
 \item $A(t;\infty,\infty,\infty)$ = $A(t-1;\ast,\ast,\ast)$,
 \item $A(t;0,\infty,\infty)$ = $A(t-1;0,\ast,\ast)$,
 \item $A(t;0,0,\ast)$ = $A(t=1;0,0,\ast)$.
 \item $A(q=s=t=1;\ast,\ast,\ast)$ = $T(3,4)$ and $A(q=s=t=1;0,\ast,\ast)$ = $P(2,-3,-2)$.
\end{enumerate}
\end{lemma}
\begin{proof}
 (1) Figure 9, Figure 14 and Figure 16 are the same.
 
 (2) By applying ambient isotopy, one can see that Figure 12, Figure 13 and Figure 15 depict the same links.
 
 (3) Replacing $t$ by $t-1$ in Figure 3 and considering Figure 11,
 one can see that $A(t;\infty,\infty,\infty)$ = $A(t-1;\ast,\ast,\ast)$.
 
 (4) Replacing $t$ by $t-1$ in Figure 6 and considering Figure 13,
 one can see that $A(t;0,\infty,\infty)$ = $A(t-1;0,\ast,\ast)$.
 
 (5) Figure 9 does not depends on $t$.
 
 (6) This comes from direct computation of the links.
 
\end{proof}
\begin{lemma}\label{lemma: 5.4}
 \begin{enumerate}
 \item det $A(t;\ast,\ast,\ast)$ = det $A(t;0,\ast,\ast)$ + det $A(t;\infty,\ast,\ast)$
 \item det $A(t;0,\ast,\ast)$ = det $A(t;0,\infty,\ast)$ + det $A(t;0,0,\ast)$  
 \item det $A(t;\infty,\ast,\ast)$ = det $A(t;\infty,0,\ast)$ + det $A(t;\infty,\infty,\ast)$  
 \item det $A(t;0,\infty,\ast)$ = det $A(t;0,\infty,0)$ + det $A(t;0,\infty,\infty)$  
 \item det $A(t;\infty,0,\ast)$ = det $A(t;0,\infty,0)$ + det $A(t;0,\infty,\infty)$  
 \item det $A(t;\infty,\infty,\ast)$ = det $A(t;0,\infty,\infty)$ + det $A(t;\infty,\infty,\infty)$  
\end{enumerate}

\end{lemma}
\begin{proof}
 Table 3.
\end{proof}
\begin{claim}\label{claim: 5.5}
 The link $A(t;\ast,\ast,\ast)$ is quasi-alternating for $s> 1$.
\end{claim}
\begin{proof}
We have shown that the link $A(t=1;\ast,\ast,\ast)$ is quasi-alternating for $s>1$. 
By induction assume $t>1$ and that the link 
$A(t-1;\ast,\ast,\ast)$ is quasi-alternating. We will show that $A(t;\ast,\ast,\ast)$ is quasi-alternating.

Since $A(t;\infty,\infty,\infty)$ = $A(t-1;\ast,\ast,\ast)$ and $A(t;0,\infty,\infty)$ = $A(t-1;0,\ast,\ast)$
are quasi-alternating by induction hypothesis, then $A(t;\infty,\infty,\ast)$ is also quasi-alternating by
Lemma \ref{lemma: 5.4}. 
Since 
$A(t;0,0,\ast)$ = $A(t=1;,0,0,\ast)$ and $A(t=1;,0,0,\ast)$ is quasi-alternating then also $A(t;0,0,\ast)$, and since
 $A(t;0,0,\ast)$ = $A(t;\infty,0,0)$ = $A(t;0,\infty,0)$ and
 $A(t;\infty,0,\infty)$ = $A(t;\infty,\infty,0)$ =\\ $A(t;0,\infty,\infty)$ by Lemma \ref{lemma: 5.3},
 then $A(t;0,\infty,\ast)$ and $A(t;\infty,0,\ast)$
 are quasi-alternating by Lemma \ref{lemma: 5.4}. Therefore, $A(t;0,\ast,\ast)$ and $A(t;\infty,\ast,\ast)$ are also 
 quasi-alternating
 by Lemma \ref{lemma: 5.4}. Finaly the link $A(t;\ast,\ast,\ast)$ is quasi-alternating by Lemma \ref{lemma: 5.4}.
 
 This completes the proof that $A(t;\ast,\ast,\ast)$ is quasi-alternating for $s>1$.
\end{proof}

\begin{claim}\label{claim: 5.6}
 $\Sigma_2(A(t;\ast,\ast,\ast))$ is an L-space for $s=1$. 
 In particular $\Sigma_2(A(t;\infty,\infty,\ast))$ is an L-space for $s=1$ and $t>1$.
\end{claim}
\begin{proof}
We have two cases
\begin{enumerate}
 \item Assume first that $t=q$. We have that if $q=s=t=1$ then $A(\ast,\ast,\ast)=T(3,4)$ and $A(0,\ast,\ast)=P(2,-3,-2)$ 
 by Lemma \ref{lemma: 5.3}, and $\Sigma_2(T(3,4))$ and $\Sigma_2(P(2,-3,-2))$ are L-spaces. By induction on $t$, assume 
 $t>1$ and $\Sigma_2(A(t-1;\ast,\ast,\ast))$ and $\Sigma_2(A(t-1;0,\ast,\ast))$ are L-spaces.
We will show that $\Sigma_2(A(t;0,\ast,\ast))$ and $\Sigma_2(A(t;\ast,\ast,\ast))$ are L-spaces.

Since $A(t;\infty,\infty,\infty)$ = $A(t-1;\ast,\ast,\ast)$ and $A(t;0,\infty,\infty)$ = $A(t-1;0,\ast,\ast)$
then $\Sigma_2(A(t;\infty,\infty,\infty))$ and $\Sigma_2(A(t;0,\infty,\infty))$ are L-spaces
by induction hypothesis. Therefore, $\Sigma_2(A(t;\infty,\infty,\ast))$ is an L-space by
Lemma \ref{lemma: 5.4} ([OSz]). 
Since 
$A(t;0,0,\ast)$ = $A(t=1;,0,0,\ast)$ and $\Sigma_2(A(t=1;,0,0,\ast))$ is an L-space, then $\Sigma_2(A(t;0,0,\ast))$ is also
an L-space. Since $A(t;0,0,\ast)$ = $A(t;\infty,0,0)$ = $A(t;0,\infty,0)$ and
 $A(t;\infty,0,\infty)$ = $A(t;\infty,\infty,0)$ = $A(t;0,\infty,\infty)$ by Lemma \ref{lemma: 5.3},
 then $\Sigma_2(A(t;0,\infty,\ast))$ and $\Sigma_2(A(t;\infty,0,\ast))$
 are L-space by Lemma \ref{lemma: 5.4}. Therefore, $\Sigma_2((A(t;0,\ast,\ast))$ and 
 $\Sigma_2(A(t;\infty,\ast,\ast))$ are also L-spaces. Finaly $\Sigma_2(A(t;\ast,\ast,\ast))$ is an L-space by 
Lemma \ref{lemma: 5.4}
 for  $t=q$.
 \item Fix $q$ and assume $t\geq q$. In case (1), we have shown that if $t=q$, then
 $\Sigma_2((A(t;0,\ast,\ast))$ and $\Sigma_2(A(t;\ast,\ast,\ast))$ are also L-spaces.
 By induction on $t$, assume 
 $t>q$ and $\Sigma_2(A(t-1;\ast,\ast,\ast))$ and $\Sigma_2(A(t-1;0,\ast,\ast))$ are L-spaces.
We will show that $\Sigma_2(A(t;0,\ast,\ast))$ and $\Sigma_2(A(t;\ast,\ast,\ast))$ are L-spaces.

Since $A(t;\infty,\infty,\infty)$ = $A(t-1;\ast,\ast,\ast)$ and $A(t;0,\infty,\infty)$ = $A(t-1;0,\ast,\ast)$
then $\Sigma_2(A(t;\infty,\infty,\infty))$ and $\Sigma_2(A(t;0,\infty,\infty))$ are L-spaces
by induction hypothesis. Therefore, $\Sigma_2(A(t;\infty,\infty,\ast))$ is an L-space by
Lemma \ref{lemma: 5.4} ([OSz]). 
Since 
$A(t;0,0,\ast)$ = $A(t=1;,0,0,\ast)$ and $\Sigma_2(A(t=1;,0,0,\ast))$ is an L-space, then $\Sigma_2(A(t;0,0,\ast))$ is also
an L-space. Since $A(t;0,0,\ast)$ = $A(t;\infty,0,0)$ = $A(t;0,\infty,0)$ and
 $A(t;\infty,0,\infty)$ = $A(t;\infty,\infty,0)$ = $A(t;0,\infty,\infty)$ by Lemma \ref{lemma: 5.3},
 then $\Sigma_2(A(t;0,\infty,\ast))$ and $\Sigma_2(A(t;\infty,0,\ast))$
 are L-space by Lemma \ref{lemma: 5.4}. Therefore, $\Sigma_2((A(t;0,\ast,\ast))$ and 
 $\Sigma_2(A(t;\infty,\ast,\ast))$ are also L-spaces by Lemma \ref{lemma: 5.4}.
Finaly $\Sigma_2(A(t;\ast,\ast,\ast))$ is an L-space by Lemma \ref{lemma: 5.4}
 for $t\geq q$.
\end{enumerate}
Since $q$ was arbitrary, then this is true for any $t$ and $q$ such that $t\geq q$.
Since $t$ and $q$ are symmetric for $A$, then $\Sigma_2(A(t;\ast,\ast,\ast))$ is an L-space for $s=1$.

\end{proof}

2) If $q>0$, $s<0$, $t>0$ or $q<0$, $s>0$, $t<0$ then the link $A(t;\ast,\ast,\ast)$ is alternating and 
therefore quasi-alternating.

3) If $q<0$, $s>0$, $t>0$ a similar argument as that of case 1) shows 
when $A(t;\ast,\ast,\ast)$ is quasi-alternating or when $\Sigma_2( A(t;\ast,\ast,\ast))$ is an L-space.

4) If $q<0$, $s<0$, $t<0$, then the link $A(t;\ast,\ast,\ast)$ is the mirror image of the link $A(t;\ast,\ast,\ast)$
when $q>0$, $s>0$, $t>0$, so a similar argument as that of case 1) shows 
when it is quasi-alternating or when $\Sigma_2( A(t;\ast,\ast,\ast))$ is an L-space.

\subsection{Proof of Theorem \ref{thm: main result1}}
The following theorem is essential in the proof of Theorem \ref{thm: main result1}.
\begin{theorem} {\rm (Theorem 3 in [MV])}\label{thm: 2.5}
 The 3-fold cyclic branched cover of $K_{[-2q,2s,-2t,2l]}$ is the 2-fold branched cover of the link
 $L(l;\ast,\ast,\ast)$.
\end{theorem}

To prove Theorem \ref{thm: main result1} we have eight cases:

(1) If $q>0$, $s>0$, $t>0$, $l>0$.
 
 \begin{figure}[H]
\centering
\def\svgwidth{0.50\columnwidth}
 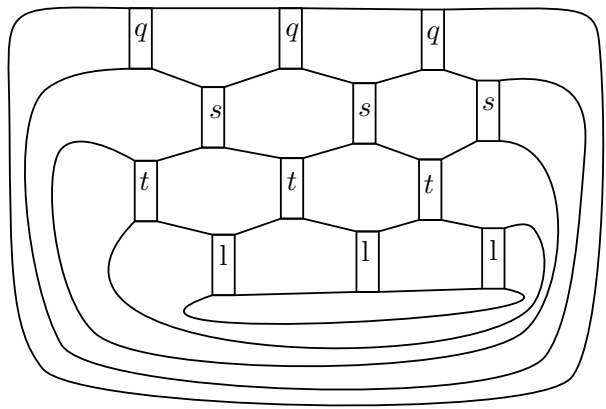
 \caption{The link $L$.}
\label{fig:The link L.}
\end{figure}

\begin{figure}[H]
\centering
\def\svgwidth{0.50\columnwidth}
 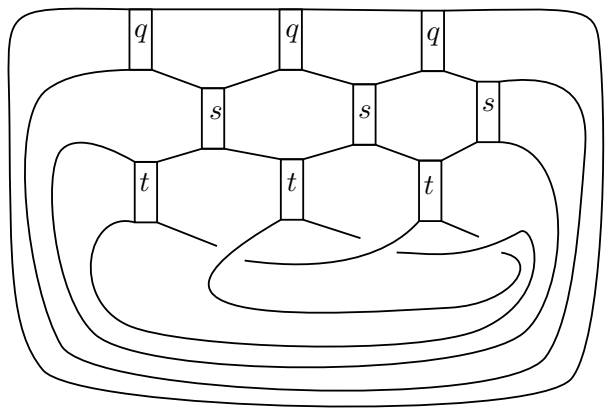
 \caption{The link $L$ for $l=1$.}
\label{fig:The link L for l=1.}
\end{figure}
\begin{figure}[H]
\centering
\def\svgwidth{0.50\columnwidth}
 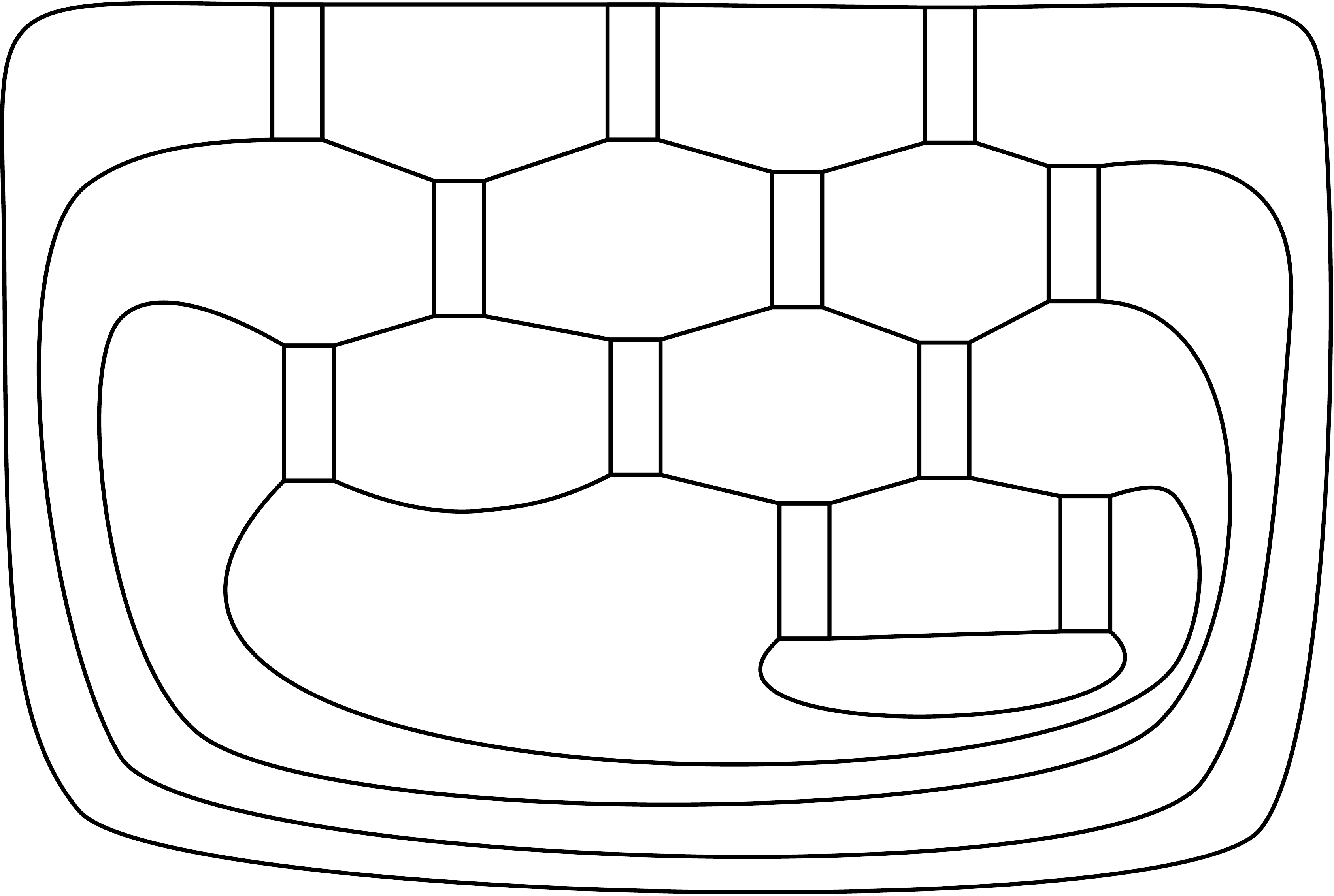
 \caption{The link $L(l;0,\ast,\ast)$.}
\label{fig:The link A.}
\end{figure}
\begin{figure}[H]
\centering
\def\svgwidth{0.50\columnwidth}
 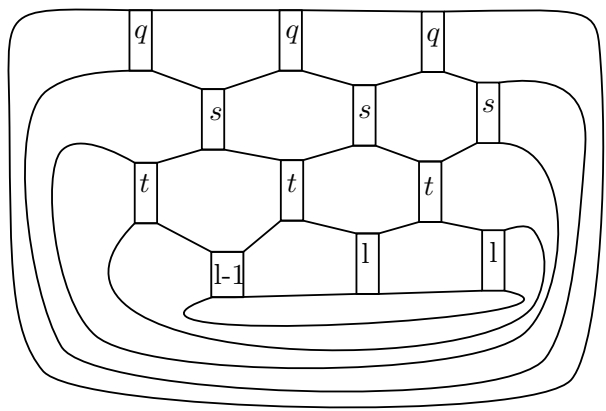
 \caption{The link $L(l;\infty,\ast,\ast)$.}
\label{fig:The link A.}
\end{figure}
\begin{figure}[H]
\centering
\def\svgwidth{0.50\columnwidth}
 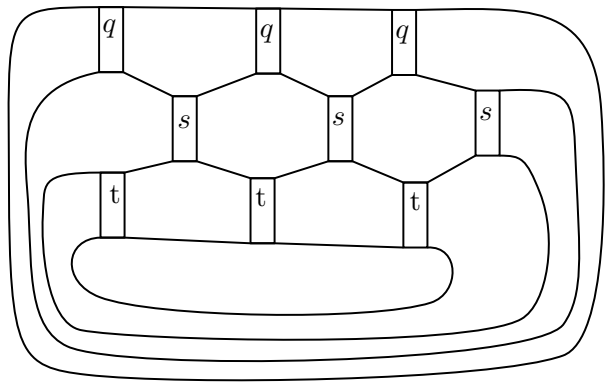
 \caption{The link $L(l;0,0,\ast)$.}
\label{fig:The link A.}
\end{figure}
\begin{figure}[H]
\centering
\def\svgwidth{0.50\columnwidth}
 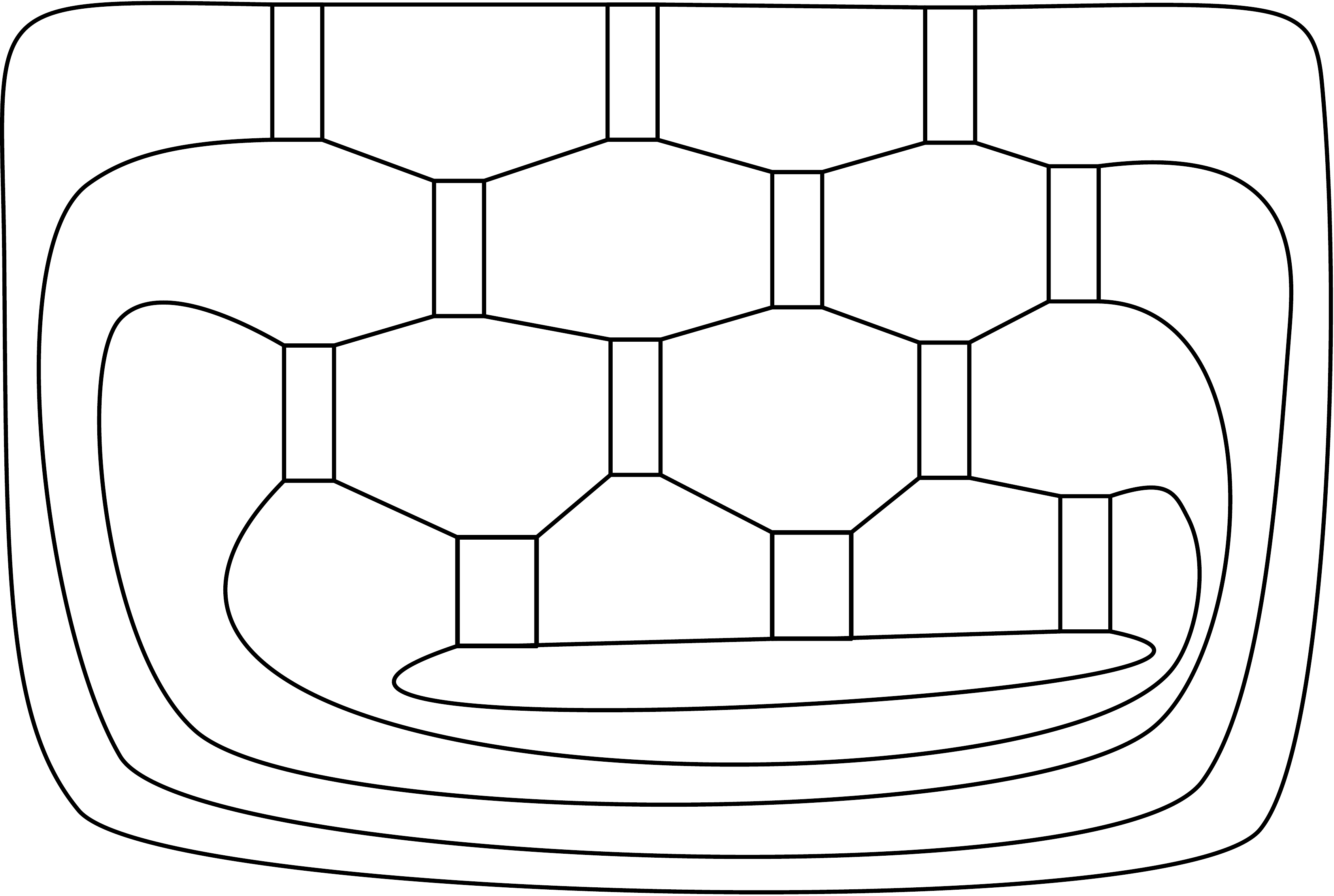
 \caption{The link $L(l;\infty,\infty,\ast)$.}
\label{fig:The link A.}
\end{figure}
\begin{figure}[H]
\centering
\def\svgwidth{0.50\columnwidth}
 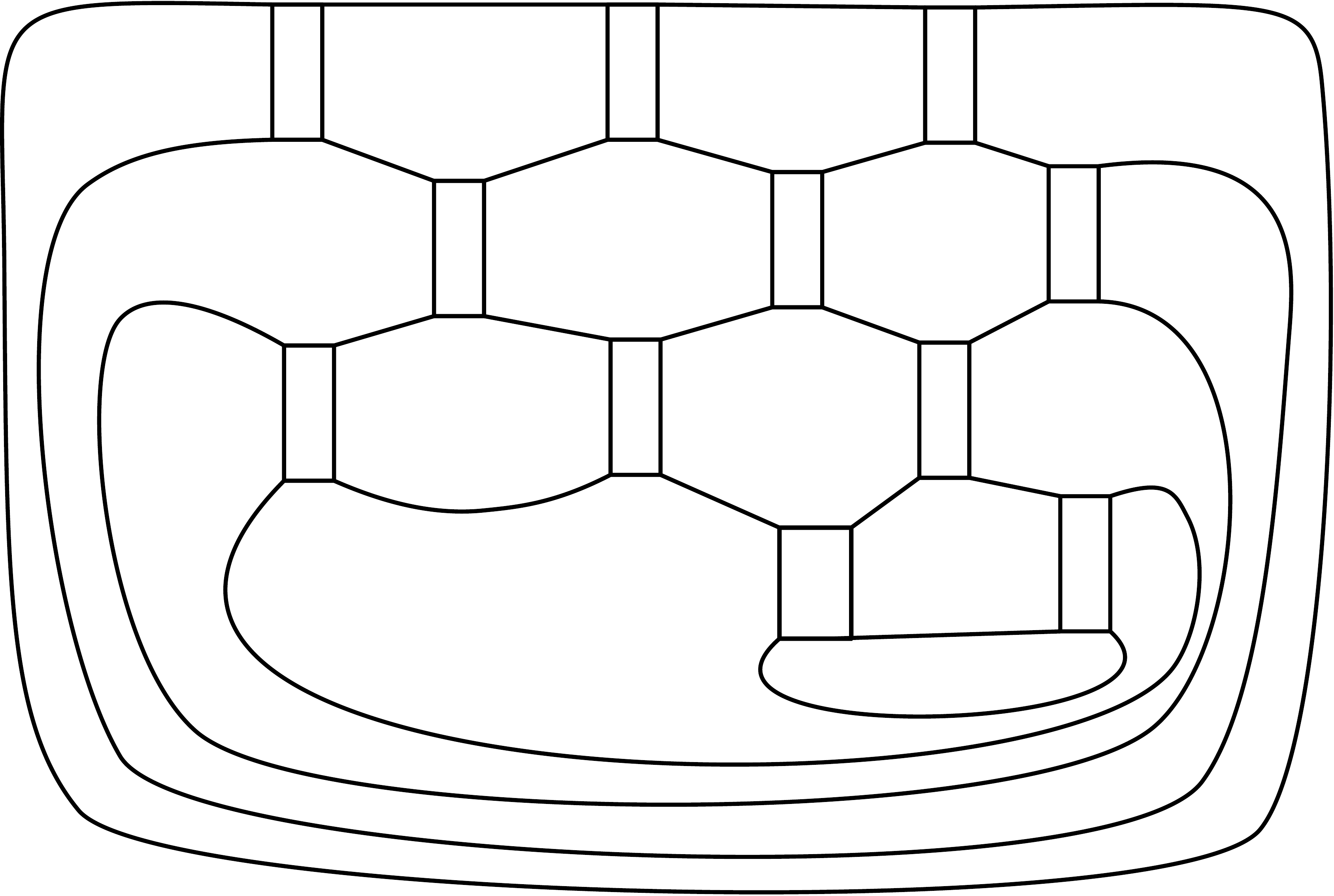
 \caption{The link $L(l;0,\infty,\ast)$.}
\label{fig:The link A.}
\end{figure}
\begin{figure}[H]
\centering
\def\svgwidth{0.50\columnwidth}
 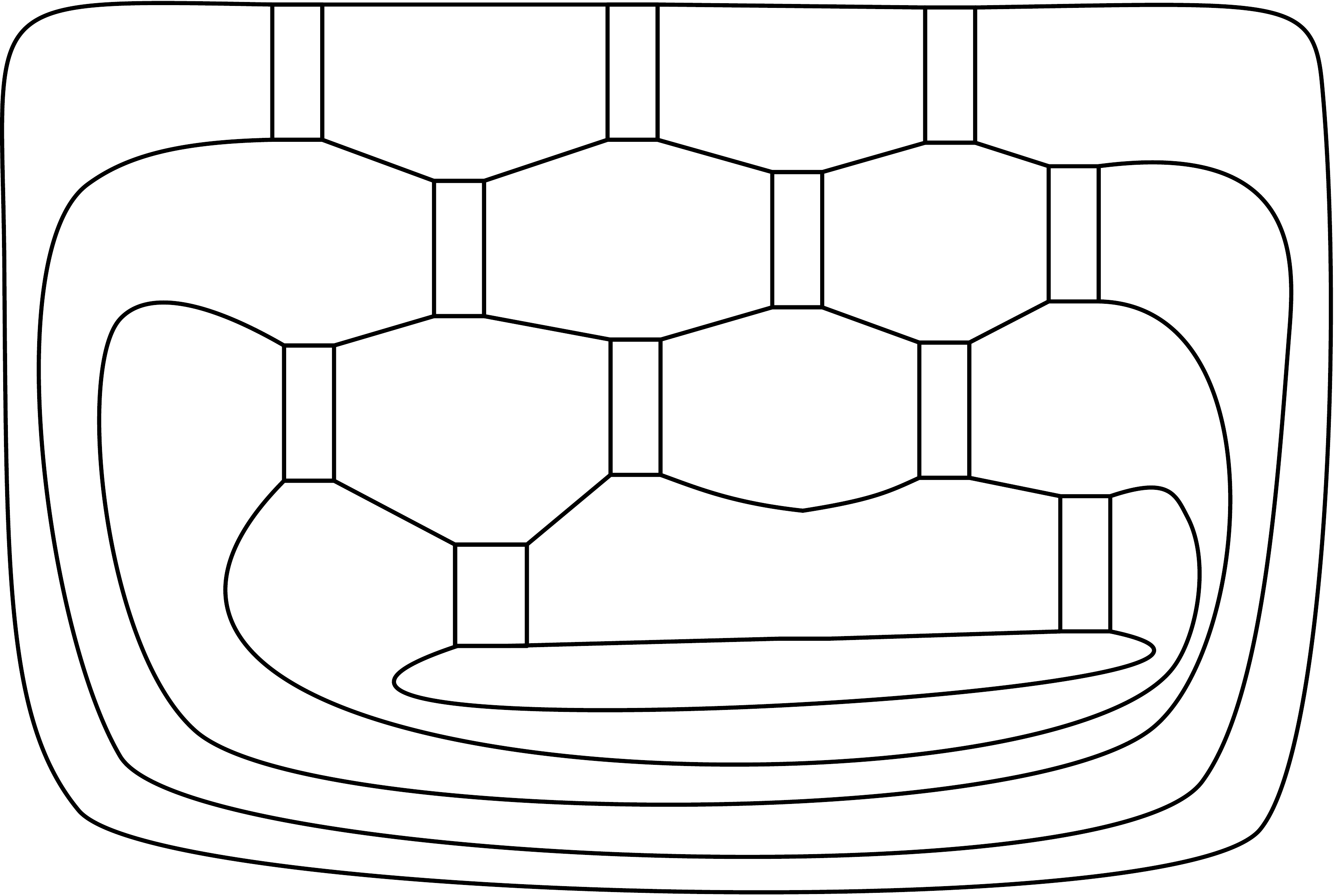
 \caption{The link $L(l;\infty,0,\ast)$.}
\label{fig:The link A.}
\end{figure}
\begin{figure}[H]
\centering
\def\svgwidth{0.50\columnwidth}
 \input{link2.eps_tex}
 \caption{The link $L(l;\infty,0,0)$.}
\label{fig:The link A.}
\end{figure}
\begin{figure}[H]
\centering
\def\svgwidth{0.50\columnwidth}
 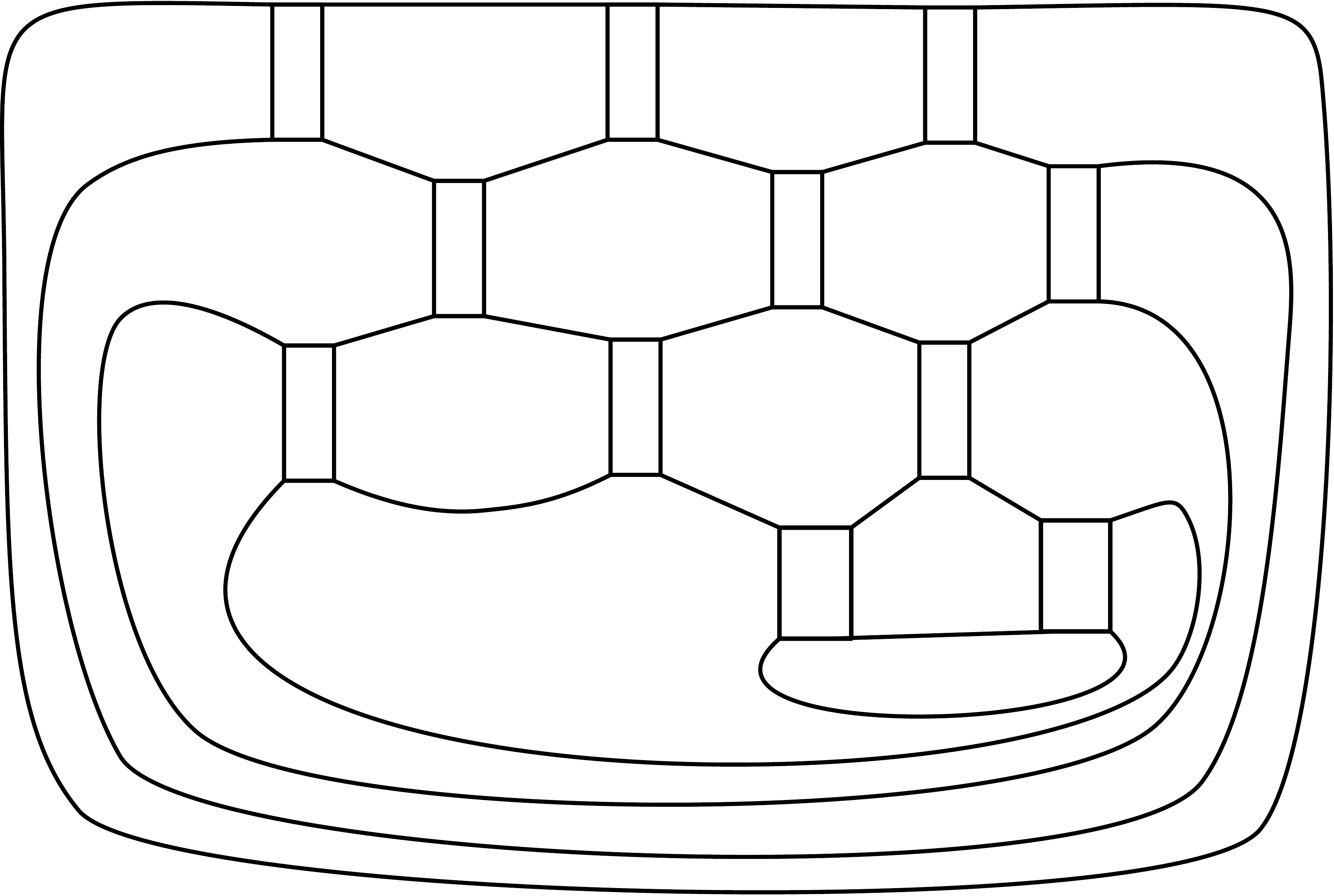
 \caption{The link $L(l;0,\infty,\infty)$.}
\label{fig:The link A.}
\end{figure}
\begin{figure}[H]
\centering
\def\svgwidth{0.50\columnwidth}
 \input{link2.eps_tex}
 \caption{The link $L(l;0,\infty,0)$.}
\label{fig:The link A.}
\end{figure}
\begin{figure}
\centering
\def\svgwidth{0.50\columnwidth}
 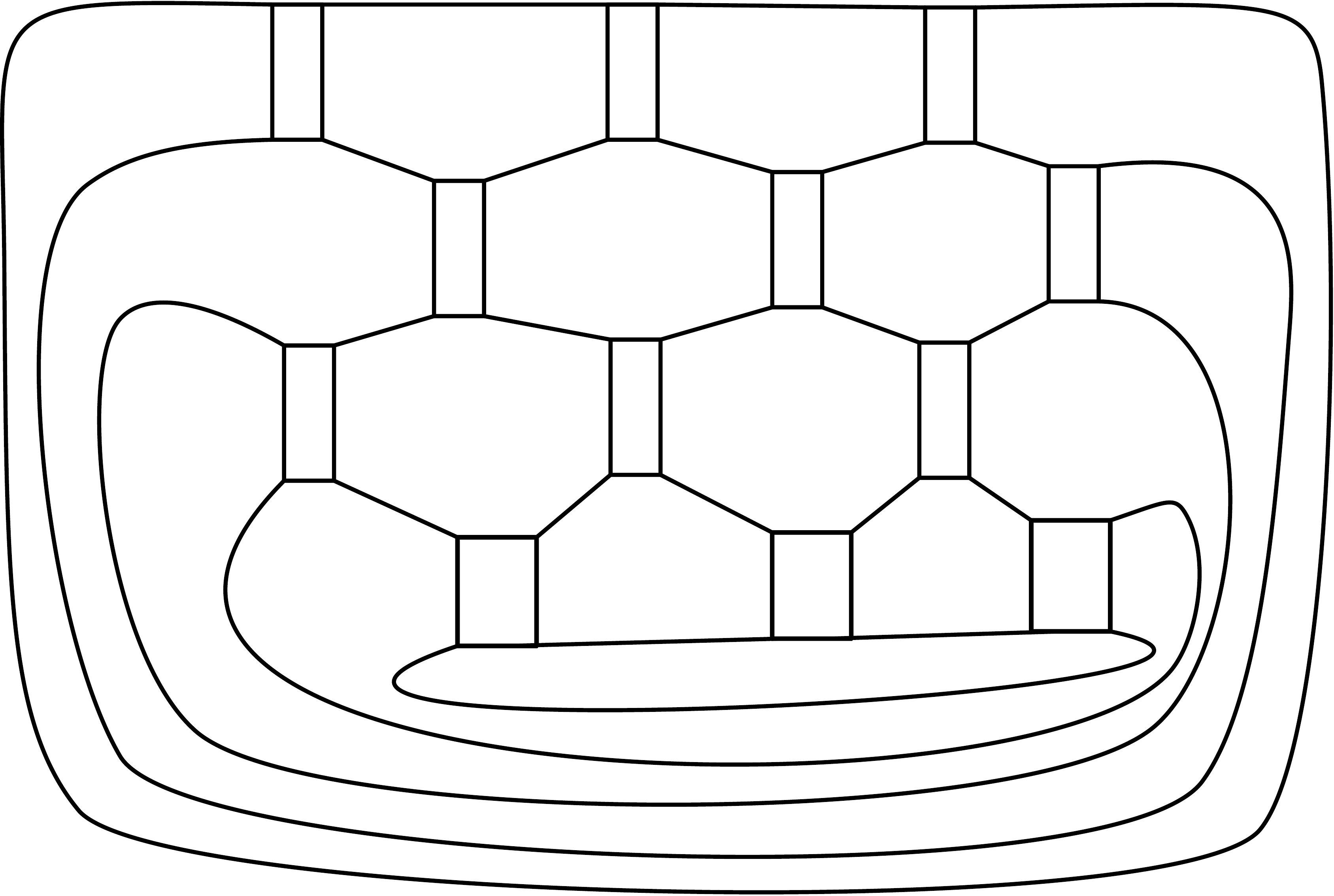
 \caption{The link $L(l;\infty,\infty,\infty)$.}
\label{fig:The link A.}
\end{figure}
\begin{figure}[H]
\centering
\def\svgwidth{0.50\columnwidth}
 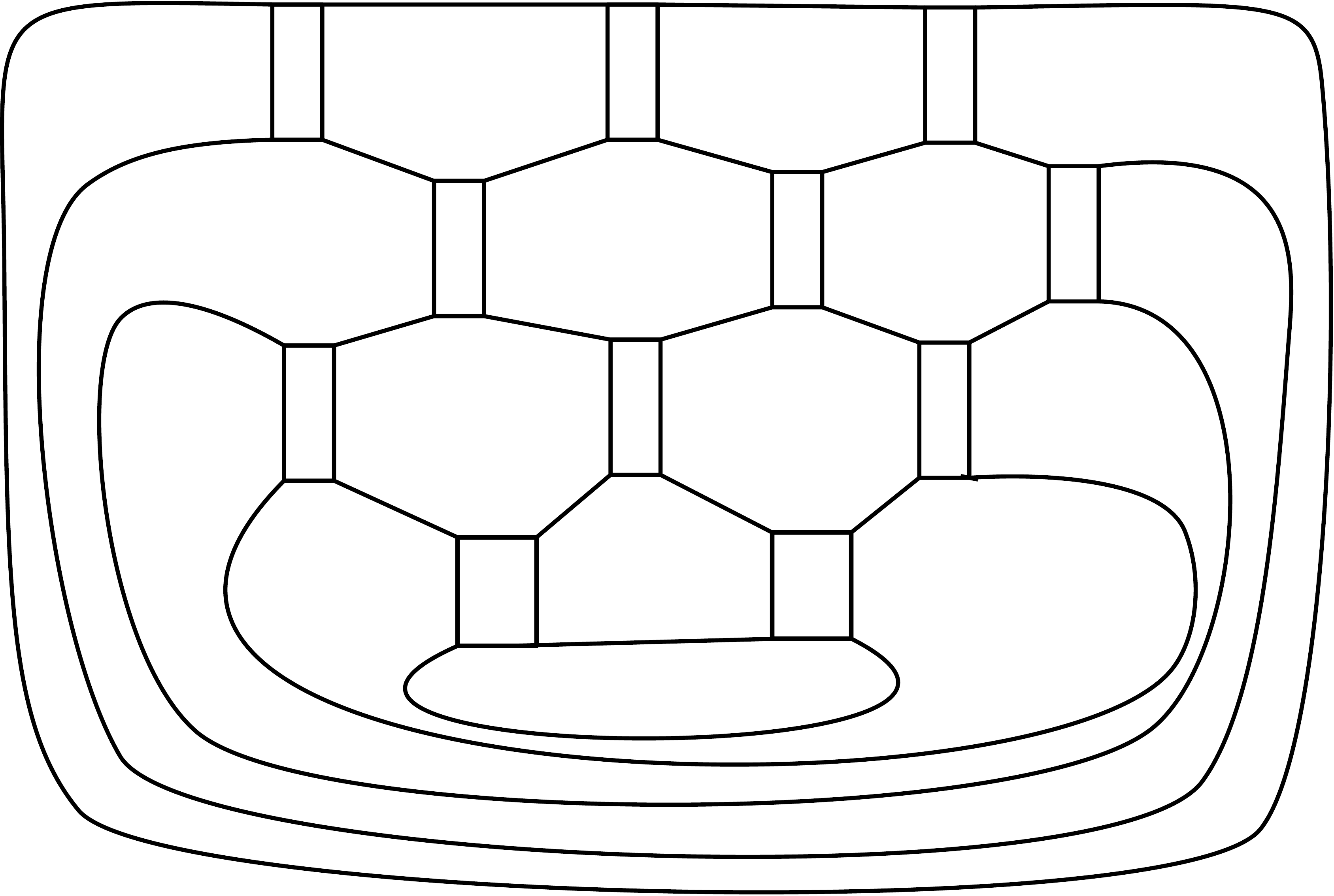
 \caption{The link $L(l;\infty,\infty,0)$.}
\label{fig:The link A.}
\end{figure}
\begin{figure}[H]
\centering
\def\svgwidth{0.50\columnwidth}
 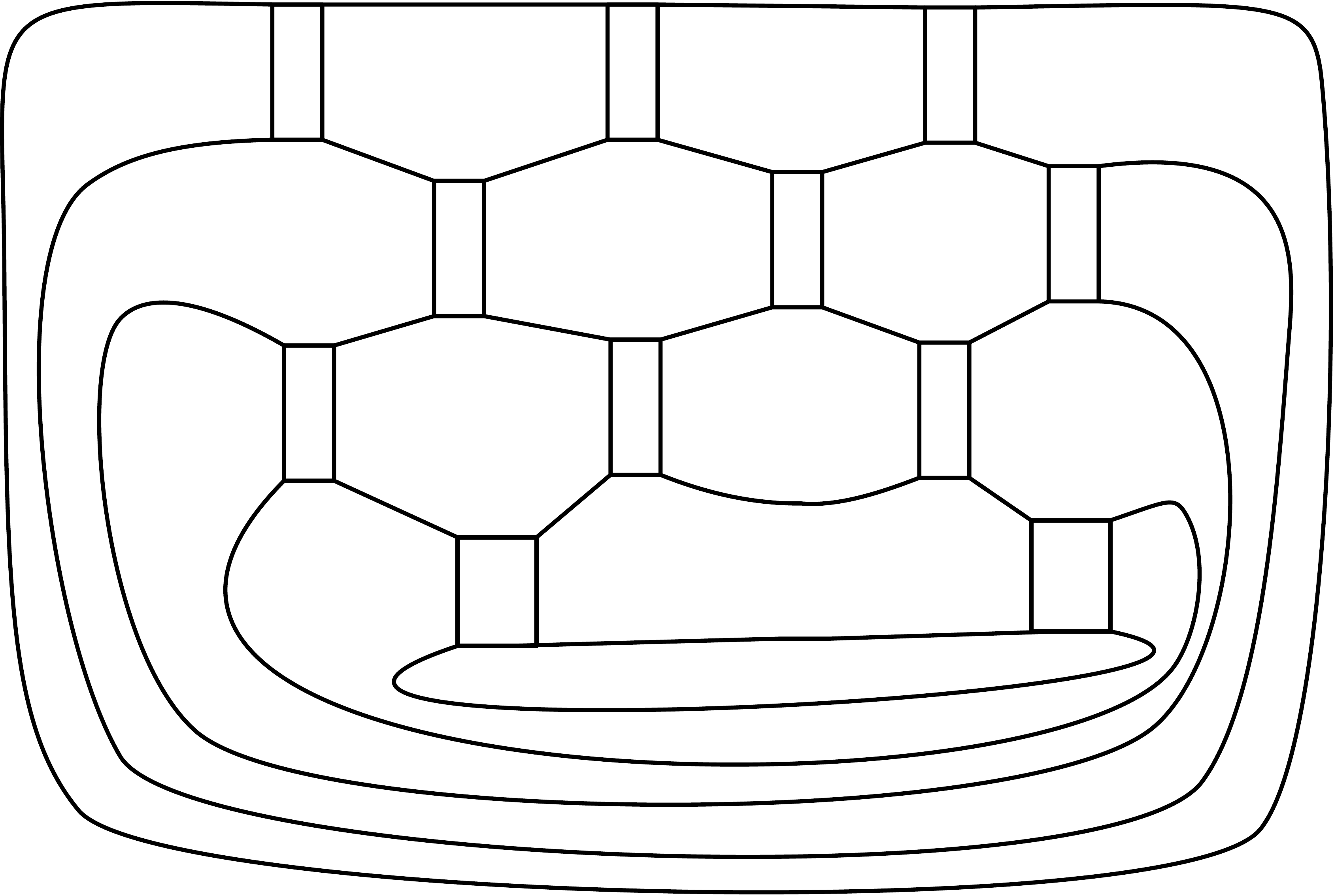
 \caption{The link $L(l;\infty,0,\infty)$.}
\label{fig:The link A.}
\end{figure}

 The Goeritz matrix of $L$ is 
 \[
  \begin{bmatrix}
    -2I & I & O & \dots  & O & \dots \\
    I & -2I & I & O & \dots  & O & \dots \\
    \vdots & \vdots & \vdots & \ddots & \vdots & \vdots & \vdots \\
    O & O \dots   &    O & I & P & I & O & \dots\\ 
     O & \dots  & O& O & I & -2I & I & O & \dots  \\
     \vdots & \vdots & \vdots & & \vdots & \ddots & \vdots & \vdots & \vdots\\
     O & \dots  & O & O & \dots O &  I & -2I & I & O \\
     O & \dots  & O & O & \dots O & O &  I & -2I & I & \\
    O & O & O & \dots  & \dots & O & O & I &Q
\end{bmatrix}
\]
where \[
   P=
  \left[ {\begin{array}{ccc}
   2s-2 & -s & -s \\
   -s & 2s-2 & -s \\
   -s & -s & 2s-2
  \end{array} } \right]
\] 
and \[
   Q=
  \left[ {\begin{array}{ccc}
   2l-1 & -l & -l \\
   -l & 2l-1 & -l \\
   -l & -l & 2l-1
  \end{array} } \right]
\]
\subsubsection{Case $l=1$}

Let $B(\ast,\ast,\ast)=L(l=1;\ast,\ast,\ast)$. 

We have the following table which gives the determinant of the corresponding link.
\begin{table}[H] \label{Table: 4}
\centering
\begin{tabular}{ | m{2cm} | m{8cm} | } 
\hline
Link & Determinant ($t>1$) \\ 
\hline
$B(\ast,\ast,\ast)$ & $(1-3q-3qs-3t+9qst)^2$ \\ 
\hline
 $B(0,\ast,\ast)$ & $2(-q-t+3qst)(1-3q-3qs-3t+9qst)$  \\ 
\hline
$B(\infty,\ast,\ast)$ & $(1-3q-3qs-3t+9qst)(1-q-3qs-t+3qst)$  \\ 
\hline
$B(0,\infty,\ast)$ & $(-q-t+3qst)(2-3q-3t-6qs+9qst)$   \\ 
\hline
 $B(0,0,\ast)$& $3(-q-t+3qst)^2$   \\ 
\hline
\end{tabular}
\caption{}
\end{table}

\begin{lemma}\label{lemma: 5.5} We have the following:
\begin{enumerate}
 \item $A(t;\ast,\ast,\ast)=B(0,0,\ast)$,
 $B(\infty,\ast,\ast)=A(t;\infty,\infty,\ast)$ and $B(0,\infty,\ast)=A(\infty,\ast,\ast)$
 \item If $q=s=t=1$ then $B(\ast,\ast,\ast)=T(3,5)$ and $B(0,\ast,\ast)=P(2,-3,-4)$.
\end{enumerate}
\end{lemma}
\begin{proof}
Recall that the link $B(\ast,\ast,\ast)$ is the following:
\begin{figure}[H]
\centering
\def\svgwidth{0.50\columnwidth}
 \input{link3.eps_tex}
 \caption{The link $B(\ast,\ast,\ast)=L(l=1,\ast,\ast,\ast)$.}
\label{fig:The link L for l=1.}
\end{figure}
Therefore, by definition of $B(0,0,\ast)$, $B(\infty,\ast,\ast)$ and $B(0,\infty,\ast)$,
the lemma follows by applying the resolution
at the right places.

\end{proof}

\begin{lemma}\label{lemma: 5.6} 
The Link  $B(\ast,\ast,\ast)$ is quasi-alternating for $s>1$, $t>1$.
\end{lemma}
\begin{proof}
Since $A(t;\ast,\ast,\ast)$ and $A(t;\infty,\infty,\ast)$ are quasi-alternating for $s> 1$ and $t>1$ 
by Claim \ref{claim: 5.5} then $B(0,0,\ast)$ and 
$B(\infty,\ast,\ast)$ are also quasi-alternating for $s> 1$ and $t>1$ by Lemma \ref{lemma: 5.5}.
By table 4, we have that  det $B(\ast,\ast,\ast)$ = det $B(0,\ast,\ast)$ + det $B(\infty,\ast,\ast)$
 and det $B(0,\ast,\ast)$ = det $B(0,\infty,\ast)$ + det $B(0,0,\ast)$ for $t> 1$, therefore $B(\ast,\ast,\ast)$ 
is quasi-alternating for $s>1$ and $t>1$.

\end{proof}
\begin{lemma}\label{lemma: 5.9}
The 2-fold branched cover of $B(\ast,\ast,\ast)$ and $B(0,\ast,\ast)$ are L-spaces for $s=1$ and $t>1$.
\end{lemma}
\begin{proof}
First, since $\Sigma_2(A(t;\ast,\ast,\ast))$ and $\Sigma_2(A(t;\infty,\infty,\ast))$ are L-spaces for $s=1$ by Claim 
\ref{claim: 5.6}, then $\Sigma_2(B(0,0,\ast))$ and 
$\Sigma_2(B(\infty,\ast,\ast))$ are also L-spaces for $s=1$ by Lemma \ref{lemma: 5.5}.
Second, we have that det $B(\ast,\ast,\ast)$ = det $B(0,\ast,\ast)$ + det $B(\infty,\ast,\ast)$
 and det $B(0,\ast,\ast)$ = det $B(0,\infty,\ast)$ + det $B(0,0,\ast)$ for $t> 1$ by Table 4, therefore 
 $\Sigma_2(B(0,\ast,\ast))$ and $\Sigma_2(B(\ast,\ast,\ast))$ 
are L-spaces for $s=1$ and $t>1$.
\end{proof}

\subsubsection{Case $l>1$}

Proceding similarly as for $A(t;\ast,\ast,\ast)$, 
we have the following table which gives the determinant of the corresponding link.
\begin{table}[H]
\centering
\begin{tabular}{ | m{3cm} | m{8cm} | } 
\hline
Link & Determinant ($l>1$)  \\ 
\hline
$L(l;\ast,\ast,\ast)$ & $(1-3ql-3qs-3lt+9lqst)^2$ \\ 
\hline
 $L(l;0,\ast,\ast)$ & $2(-t-q+3qst)(1-3lq-3qs-3lt+9lqst)$  \\ 
\hline
$L(l;\infty,\ast,\ast)$ & $(1-3lq-3qs-3lt+9lqst)(1+2q-3lq-3qs+2t-3lt-6qst+9lqst)$  \\ 
\hline
$L(l;0,\infty,\ast)$ & $(-q-t+3qst)(2+3q-6lq-6qs+3t-6lt-9qst+18lqst)$   \\ 
\hline
$L(l;\infty,0,\ast)$ & $(-q-t+3qst)(2+3q-6lq-6qs+3t-6lt-9qst+18lqst)$  \\ 
\hline
 $L(l;0,0,\ast)$& $3(-q-t+3qst)^2$   \\ 
\hline
$L(l;\infty,\infty,\ast)$ & $(1+3q-3lq-3qs+3t-3lt-9qst+9lqst)(1+q-3lq-3qs+t-3lt-3qst+9lqst)$  \\ 
\hline
$L(l;\infty,\infty,\infty)$ & $(1+3q-3lq-3qs+3t-3lt-9qst+9lqst)^2$  \\ 
\hline
$L(l;0,\infty,\infty)$ & $2(-q-t+3qst)(1+3q-3lq-3qs+3t-3lt-9qst+9lqst)$  \\ 
\hline
\end{tabular}
\caption{}
\end{table}

\begin{lemma}\label{lemma: 5.7}
\begin{enumerate}
 \item $L(l;0,0,\ast)$ = $L(l;\infty,0,0)$ = $L(l;0,\infty,0)$, 
 \item $L(l;\infty,0,\infty)$ = $L(l;\infty,\infty,0)$ = $L(l;0,\infty,\infty)$, 
 \item $L(l;\infty,\infty,\infty)$ = $L(l-1;\ast,\ast,\ast)$,
 \item $L(l;0,\infty,\infty)$ = $L(l-1;0,\ast,\ast)$,
 \item $L(l;0,0,\ast)$ = $L(l=1;,0,0,\ast)$.
\end{enumerate}

\end{lemma}
\begin{proof}
(1) Figure 21, Figure 25 and Figure 27 are the same.
 
 (2) By applying ambient isotopy, one can see that Figure 26, Figure 29 and Figure 30 depict the same links.
 
 (3) Replacing $l$ by $l-1$ in Figure 17 and considering Figure 28, 
 one can see that $L(l;\infty,\infty,\infty)$ = $L(l-1;\ast,\ast,\ast)$.
 
 (4) Replacing $l$ by $l-1$ in Figure 19 and considering Figure 26, 
 one can see that $L(l;0,\infty,\infty)$ = $L(l-1;0,\ast,\ast)$.
 
 (5) Figure 21 does not depends on $l$.

\end{proof}
\begin{lemma}\label{lemma: 5.8}
 \begin{enumerate}
 \item det $L(l;\ast,\ast,\ast)$ = det $L(l;0,\ast,\ast)$ + det $L(l;\infty,\ast,\ast)$
 \item det $L(l;0,\ast,\ast)$ = det $L(l;0,\infty,\ast)$ + det $L(l;0,0,\ast)$  
 \item det $L(l;\infty,\ast,\ast)$ = det $L(l;\infty,0,\ast)$ + det $L(l;\infty,\infty,\ast)$  
 \item det $L(l;0,\infty,\ast)$ = det $L(l;0,\infty,0)$ + det $L(l;0,\infty,\infty)$  
 \item det $L(l;\infty,0,\ast)$ = det $L(l;0,\infty,0)$ + det $L(l;0,\infty,\infty)$  
 \item det $L(l;\infty,\infty,\ast)$ = det $L(l;0,\infty,\infty)$ + det $L(l;\infty,\infty,\infty)$  
\end{enumerate}

\end{lemma}
\begin{proof}
 See Table 5.
\end{proof}
\begin{proof}[Proof of Theorem \ref {thm: main result1}]
The proof of this theorem will be split into 3 claims.
\begin{claim}\label{claim: 5.12}
If $s>1$ and $t>1$, then the 2-fold branched cover of $L(l;\ast,\ast,\ast)$ is an L-space.
\end{claim}
\begin{proof}
We have shown that the link $L(l=1;\ast,\ast,\ast)=B(\ast,\ast,\ast)$ is quasi-alternating for $s>1$ and $t>1$.
By induction on $l$, assume $l>1$ and
the link 
$L(l-1;\ast,\ast,\ast)$ is quasi-alternating for $s>1$ and $t>1$.
We will show that $L(l;\ast,\ast,\ast)$ is quasi-alternating for $s>1$ and $t>1$.

Since $L(l;\infty,\infty,\infty)$ = $L(l-1;\ast,\ast,\ast)$ and $L(l;0,\infty,\infty)$ = $L(l-1;0,\ast,\ast)$
are quasi-alternating by induction hypothesis, then $L(l;\infty,\infty,\ast)$ is also quasi-alternating by 
Lemma \ref{lemma: 5.8}. Since 
$L(l;0,0,\ast)$ = $L(l=1;,0,0,\ast)$, and $L(l=1;,0,0,\ast)$ is quasi-alternating, then $L(l;0,0,\ast)$ is also 
quasi-alternating. We have that
 $L(l;0,0,\ast)$ = $L(l;\infty,0,0)$ = $L(l;0,\infty,0)$ and
 $L(l;\infty,0,\infty)$ = $L(l;\infty,\infty,0)$ = $L(l;0,\infty,\infty)$ by Lemma \ref{lemma: 5.7},
 then $L(l;0,\infty,\ast)$ and $L(l;\infty,0,\ast)$ are quasi-alternating by 
 Lemma \ref{lemma: 5.8}. Therefore, $L(l;0,\ast,\ast)$ and $L(l;\infty,\ast,\ast)$ are also quasi-alternating
 by Lemma \ref{lemma: 5.8}. Finaly the link $L(l;\ast,\ast,\ast)$ is quasi-alternating by Lemma \ref{lemma: 5.8}
 for $s>1$ and $t>1$.
 
Therefore, the 2-fold branched cover of $L(l;\ast,\ast,\ast)$ is an L-space for $s>1$ and $t>1$ by
Theorem \ref{thm: 2.4}.
\end{proof}
\begin{claim}
 If $s=t=1$, then the 2-fold branched cover of $L(l;\ast,\ast,\ast)$ is an L-space.
\end{claim}
\begin{proof}
We have two cases
\begin{enumerate}
 \item Assume first that $l=q$. We have that if $q=s=t=l=1$ then $B(\ast,\ast,\ast)=T(3,5)$ and $B(0,\ast,\ast)=P(2,-3,-4)$ 
 by Lemma \ref{lemma: 5.5}, and $\Sigma_2(T(3,5))$ and $\Sigma_2(P(2,-3,-4))$ are L-spaces. By induction on $l$, assume 
 $l>1$ and $\Sigma_2(L(l-1;\ast,\ast,\ast))$ and $\Sigma_2(L(l-1;0,\ast,\ast))$ are L-spaces.
We will show that $\Sigma_2(L(l;0,\ast,\ast))$ and $\Sigma_2(L(l;\ast,\ast,\ast))$ are L-spaces.

Since $L(l;\infty,\infty,\infty)$ = $L(l-1;\ast,\ast,\ast)$ and $L(l;0,\infty,\infty)$ = $L(l-1;0,\ast,\ast)$
then \\$\Sigma_2(L(l;\infty,\infty,\infty))$ and $\Sigma_2(L(l;0,\infty,\infty))$ are L-spaces
by induction hypothesis. Therefore $\Sigma_2(L(l;\infty,\infty,\ast))$ is an L-space by
Lemma \ref{lemma: 5.8} ([OSz]). 
Since 
$L(l;0,0,\ast)$ = $L(l=1;,0,0,\ast)$ and $\Sigma_2(L(l=1;,0,0,\ast))$ is an L-space, then $\Sigma_2(L(l;0,0,\ast))$ is also
an L-space. Since $L(l;0,0,\ast)$ = $L(l;\infty,0,0)$ = $L(l;0,\infty,0)$ and
 $L(l;\infty,0,\infty)$ = $L(l;\infty,\infty,0)$ = $L(l;0,\infty,\infty)$ by Lemma \ref{lemma: 5.7},
 then $\Sigma_2(L(l;0,\infty,\ast))$ and $\Sigma_2(L(l;\infty,0,\ast))$
 are L-space by Lemma \ref{lemma: 5.8}. Therefore, $\Sigma_2((L(l;0,\ast,\ast))$ and 
 $\Sigma_2(L(l;\infty,\ast,\ast))$ are also L-spaces by Lemma \ref{lemma: 5.8}.
 Finaly $\Sigma_2(L(l;\ast,\ast,\ast))$ is an L-space by Lemma \ref{lemma: 5.8}.
 \item Fix $q$ and assume $l\geq q$. In case (1), we have shown that if $l=q$, then
 $\Sigma_2((L(l;0,\ast,\ast))$ and $\Sigma_2(L(l;\ast,\ast,\ast))$ are also L-spaces.
 By induction on $l$, assume 
 $l>q$ and $\Sigma_2(L(l-1;\ast,\ast,\ast))$ and $\Sigma_2(L(l-1;0,\ast,\ast))$ are L-spaces.
We will show that $\Sigma_2(L(l;0,\ast,\ast))$ and $\Sigma_2(L(l;\ast,\ast,\ast))$ are L-spaces.

Since $L(l;\infty,\infty,\infty)$ = $L(l-1;\ast,\ast,\ast)$ and $L(l;0,\infty,\infty)$ = $L(l-1;0,\ast,\ast)$
then\\ $\Sigma_2(L(l;\infty,\infty,\infty))$ and $\Sigma_2(L(l;0,\infty,\infty))$ are L-spaces
by induction hypothesis. Therefore $\Sigma_2(L(l;\infty,\infty,\ast))$ is an L-space by
Lemma \ref{lemma: 5.8} ([OSz]). 
Since 
$L(l;0,0,\ast)$ = $L(l=1;,0,0,\ast)$ and $\Sigma_2(L(l=1;,0,0,\ast))$ is an L-space, then $\Sigma_2(L(l;0,0,\ast))$ is also
an L-space. Since $L(l;0,0,\ast)$ = $L(l;\infty,0,0)$ = $L(l;0,\infty,0)$ and
 $L(l;\infty,0,\infty)$ = $L(l;\infty,\infty,0)$ = $L(l;0,\infty,\infty)$ by Lemma \ref{lemma: 5.7},
 then $\Sigma_2(L(l;0,\infty,\ast))$ and $\Sigma_2(L(l;\infty,0,\ast))$
 are L-spaces by Lemma \ref{lemma: 5.8}. Therefore, $\Sigma_2((L(l;0,\ast,\ast))$ and 
 $\Sigma_2(L(l;\infty,\ast,\ast))$ are also L-spaces by Lemma \ref{lemma: 5.8}.
 Finaly $\Sigma_2(L(l;\ast,\ast,\ast))$ is an L-space by Lemma \ref{lemma: 5.8} for $l\geq q$.
\end{enumerate}
Since $q$ was arbitrary, then this is true for any $l$ and $q$ such that $l\geq q$.
Since $l$ and $q$ are symmetric for $L$, then $\Sigma_2(L(l;\ast,\ast,\ast))$ is an L-space for $s=t=1$.
\end{proof}
\begin{claim}
 If $s=1$ and $t>1$, then $\Sigma_2(L(l;\ast,\ast,\ast))$ is an L-space.
\end{claim}
\begin{proof}
 Similarly as the proofs of the two first claims, we will proceed by induction. 
 We have shown that if $s=1$ and $t>1$, 
 then $\Sigma_2(L(l=1;\ast,\ast,\ast))=\Sigma_2(B(\ast,\ast,\ast))$ and $\Sigma_2(L(l=1;0,\ast,\ast))=\Sigma_2(B(0,\ast,\ast))$ 
 are L-spaces by Lemma \ref{lemma: 5.9}. By induction on $l$, assume 
 $l>1$ and $\Sigma_2(L(l-1;\ast,\ast,\ast))$ and $\Sigma_2(L(l-1;0,\ast,\ast))$ are L-spaces.
We will show that $\Sigma_2(L(l;0,\ast,\ast))$ and $\Sigma_2(L(l;\ast,\ast,\ast))$ are L-spaces.

Since $L(l;\infty,\infty,\infty)$ = $L(l-1;\ast,\ast,\ast)$ and $L(l;0,\infty,\infty)$ = $L(l-1;0,\ast,\ast)$
then $\Sigma_2(L(l;\infty,\infty,\infty))$ and $\Sigma_2(L(l;0,\infty,\infty))$ are L-spaces
by induction hypothesis. Therefore $\Sigma_2(L(l;\infty,\infty,\ast))$ is an L-space by
Lemma \ref{lemma: 5.8} ([OSz]). 
Since 
$L(l;0,0,\ast)$ = $L(l=1;,0,0,\ast)$ and $\Sigma_2(L(l=1;,0,0,\ast))$ is an L-space, then $\Sigma_2(L(l;0,0,\ast))$ is also
an L-space. Since $L(l;0,0,\ast)$ = $L(l;\infty,0,0)$ = $L(l;0,\infty,0)$ and
 $L(l;\infty,0,\infty)$ = $L(l;\infty,\infty,0)$ = $L(l;0,\infty,\infty)$ by Lemma \ref{lemma: 5.7},
 then \\$\Sigma_2(L(l;0,\infty,\ast))$ and $\Sigma_2(L(l;\infty,0,\ast))$
 are L-space by Lemma \ref{lemma: 5.8}. Therefore, $\Sigma_2((L(l;0,\ast,\ast))$ and 
 $\Sigma_2(L(l;\infty,\ast,\ast))$ are also L-spaces by Lemma \ref{lemma: 5.8}. 
 Finaly $\Sigma_2(L(l;\ast,\ast,\ast))$ is an L-space by Lemma \ref{lemma: 5.8}.
\end{proof}
Since $s$ and $t$ are symmetric for $L(l;\ast,\ast,\ast)$, then $\Sigma_2(L(l;\ast,\ast,\ast))$ is also an L-space for $s>1$ 
and $t=1$.

This complete the proof of Theorem \ref{thm: main result1} by Theorem \ref{thm: 2.5} for $q>0$, $s>0$, $t>0$ and $l>0$.
\end{proof}

(2) If $q<0$, $s>0$, $t<0$ and $l>0$ then the link $L$ is alternating, and therefore quasi-alternating, so 
Theorem \ref{thm: main result1} is a consequence of Theorem \ref{thm: 2.5}.

(3) If $q>0$, $s<0$, $t>0$ and $l>0$ then the link $A$ is alternating and therefore quasi-alternating. Hence the proof of
Theorem \ref{thm: main result1} is dealt with as in case (1).

(4) If $q<0$, $s<0$, $t<0$ and $l>0$ then the link $A$ is the mirror image of the link $A$ in the case (1), so a similar proof
to that used in case (1) can be used to prove this case.

(5) If $q>0$, $s>0$, $t<0$ and $l>0$ then as $q$ and $l$ are symmetric for $L(l;\ast,\ast,\ast)$, 
and $s$ and $t$ are also symmetric for $L(l;\ast,\ast,\ast)$, this case of Theorem \ref{thm: main result1}
follows as in case (3).

(6) If $q>0$, $s<0$, $t<0$ and $l<0$ then as $q$ and $l$ are symmetric for $L(l;\ast,\ast,\ast)$,
and $s$ and $t$ are also symmetric for $L(l;\ast,\ast,\ast)$,
this case of Theorem \ref{thm: main result1}
follows as in case (4).

For (7), (8) the proof of Theorem \ref{thm: main result1} is done in a similar way as case (1).

\begin{remark}
{\rm 
\begin{enumerate}
 \item In general, it is not true that for every two-bridge knot $K$ the 3-fold cyclic branched cover is an $L$-space. 
Boileau-Boyer-Gordon [BBG] have shown that the 3-fold cyclic branched cover of some familly of strongly quasipositive
2-bridge knots is not an $L$-space. They also show that for some familly of genus 2 strongly quasipositive 2-bridge 
knots the 3-fold cyclic branched cover is an $L$-space.
\item Claim \ref{claim: 5.12} gives an infinite family of quasi-alternating links whose 2-fold branched cover has 
none left-orderable fundamental group.
\end{enumerate}

}
\end{remark}


\begin{thebibliography}{alpha}
\bibitem[BBG]{BBG} 
  Michel Boileau, Steven Boyer, Cameron McA. Gordon,\textit{ Branched covers of quasipositive links and L-spaces}, arxiv. 

\bibitem[BGW]{BGW}
Steven Boyer, Cameron McA. Gordon, and Liam Watson
 \textit{On L-spaces and left-orderable fundamental
groups}, Math. Ann. {\bf 356} (2013), 1213--1245.
 
 \bibitem[BRW]{BRW} 
Steven Boyer, Dale Rolfsen, and Bert Wiest,
\textit{Orderable 3-manifold groups}, Ann. Inst. Fourier  
{\bf 55} (2005), 243--288. 

\bibitem[BZ]{BZ} 
 Steven Boyer and Xingru Zhang,
\textit{Cyclic surgery and boundary slopes}, {\bf Geometric topology}, W.
Kazez ed., AMS/IP Studies in Advanced Mathematics 2 (1996), 62--79.

 \bibitem[CR]{CR} 
 Adam Clay and Dale Rolfsen, 
\textit{Ordered groups and topology}, 2015.

 \bibitem[DPT]{DPT} 
 M. Düabkowski, J. Przytycki, and A. Togha,
\textit{Non-left-orderable 3-manifold groups}, Canadian
Math. Bull., {\bf 48(1)} 32--40, 2005. 

 \bibitem[GL]{GL} 
 Cameron McA. Gordon, Tye Lidman,
\textit{Taut foliations, left-orderability, and cyclic branched covers}, Acta Math. Vietnam {\bf 39} (2014), no.4, 599--635. 

 \bibitem[Hu]{Hu} 
Ying Hu, 
\textit{The left-orderability and the cyclic branched coverings}, Algebraic and Geometric Topology {\bf 15} (2015) 399--413. 
\bibitem[K]{K} 
 Taizo Kanenobu, 
\textit{Genus and Kaufman Polynomial of a 2-bridge knot}, Osaka J. Math., {\bf 29} (1992), 635--651.

 \bibitem[Ka]{Ka} 
 A. Kawauchi, 
\textit{A survey of knot theory}, Birkhauser Verlag, Basel, 1996.

\bibitem[L]{L} 
 W.B. Raymond Lickorish, 
\textit{An Introduction to Knot Theory}, Graduate texts in mathematics, Springer, 1997.

 
\bibitem[MV]{MV} 
 Michele Mulazzani and Andrei Vesnin,
\textit{Generalized Takahashi manifolds}, Osaka J. Math.  
{\bf 39} (2002), 705--721. 

\bibitem[OSz]{OSz} 
 Peter Ozsvath and Zoltan Szabo,
\textit{On the Heegaard Floer homolgy of branched double-covers}, Adv. Math.  
{\bf 194} (2005), 1-33. 
\bibitem[OSz06]{OSz} 
 Peter Ozsvath and Zoltan Szabo, \textit{Introduction to Heegaard Floer theory}, Clay Mathematics Institute.
 {\bf 5} (2006), 3--28. 
\bibitem[P]{P} 
 T. Peters, 
\textit{On L-spaces and non left-orderable 3-manifold groups}, preprint, arXiv.

 \bibitem[Pr]{Pr} 
 Jozef H. Przytycki, 
\textit{From Goeritz matrix to quasi-alternating links}, preprint, arXiv.

\bibitem[R]{R} 
 Dale Rolfsen, 
\textit{knots and Links}, Publish or Perish, 1976.

\bibitem[Te]{Te} 
 Masakazu Teragaito,
\textit{Cyclic branched covers of alternating knots and L-spaces}, Bull. Korean Math. Soc.  
{\bf 52} (2015), 1139--1148. 

 \bibitem[Tra]{Tra} 
 Anh T. Tran, 
\textit{On left-orderability and cyclic branched coverings}, Journal of the 
Mathematical Society of Japan {\bf 67} (2015), no. 3, 1169--1178. 

 \end{thebibliography}
\end{document}